\documentclass[11pt,a4paper]{article}

\usepackage[a4paper,margin=20mm]{geometry}
\usepackage[T1]{fontenc}
\usepackage{lmodern}
\usepackage{microtype}

\usepackage{authblk}

\usepackage[round,authoryear]{natbib}

\usepackage{graphicx}
\usepackage{booktabs}
\usepackage{multirow}
\usepackage{subcaption}
\usepackage{verbatim}

\usepackage{tikz}
\usetikzlibrary{arrows.meta,positioning,calc,petri}

\usepackage{amssymb}
\usepackage{amsmath}
\usepackage{amsthm}
\newtheorem{proposition}{Proposition}
\newcommand{\restatedpropnum}{}
\newtheorem*{restatedpropinner}{Proposition~\restatedpropnum}
\newenvironment{restatedprop}[1]
  {\renewcommand{\restatedpropnum}{#1}\begin{restatedpropinner}}
  {\end{restatedpropinner}}

\usepackage{enumitem}
\setlist[itemize,enumerate]{
    topsep=3pt,
    partopsep=0pt,
    itemsep=0pt,
    parsep=0pt
}

\AtBeginDocument{
    \setlength{\abovedisplayskip}{8pt plus 2pt minus 2pt}
    \setlength{\belowdisplayskip}{8pt plus 2pt minus 2pt}
    \setlength{\abovedisplayshortskip}{4pt plus 2pt}
    \setlength{\belowdisplayshortskip}{5pt plus 2pt minus 2pt}
}

\usepackage[dvipsnames]{xcolor}

\newcommand{\red}[1]{#1} 

\newcommand{\E}{\mathcal{E}}
\newcommand{\V}{\mathcal{V}}
\newcommand{\R}{\mathcal{R}}
\newcommand{\T}{\mathcal{T}}
\newcommand{\J}{\mathcal{J}}

\newcommand{\Sset}{\mathcal{S}}

\newcommand{\dres}[1]{\tilde{d}_{#1}}  
\newcommand{\feval}{\hat f}  
\newcommand{\Pred}{\operatorname{Pred}}
\newcommand{\Suc}{\operatorname{Suc}}

\usepackage[hidelinks]{hyperref}

\title{Two Communities, One Problem: Branch-and-Bound vs.\ \texorpdfstring{$A^*$}{A*} for Resource-Constrained Project Scheduling}

\author[1]{Izack Cohen}
\author[2]{Dominik Kress}
\author[2]{Lara Nehrke}
\author[3]{Erwin Pesch}
\author[1]{Shimrit Shtern}

\affil[1]{Faculty of Engineering, Bar-Ilan University, Max and Anna Webb Street, Ramat Gan 5290002, Israel}
\affil[2]{Business Administration, especially Procurement and Production, Helmut Schmidt University/University of the Federal Armed Forces Hamburg, Friedrich-Ebert-Damm 245, D-22159 Hamburg, Germany}
\affil[3]{Management Information Science, University of Siegen, Kohlbettstr.~15, D-57072 Siegen, Germany}
\affil[ ]{\texttt{izack.cohen@biu.ac.il}, \texttt{dominik.kress@hsu-hh.de}, \texttt{lara.nehrke@hsu-hh.de}, \texttt{erwin.pesch@uni-siegen.de}, \texttt{shimrit.shtern@biu.ac.il}}

\date{}

\begin{document}

\maketitle

\begin{abstract}
Branch-and-Bound (B\&B) and $A^*$ are popular exact search methods originating from Operations Research and Artificial Intelligence, respectively. Although the two methods share fundamental search principles, their different ways of exploring the search space can lead to substantially different computational performance. We investigate these differences using the Resource-Constrained Project Scheduling Problem (RCPSP) as an exemplary problem setting. We first establish the relationships between mixed-integer linear programming, state-action, and timed Petri-net representations. Since B\&B and $A^*$ admit numerous algorithmic enhancements, we focus on relatively basic implementations. This allows us to examine the fundamental differences between the methods rather than compare the performance of state-of-the-art solvers. By using the same lower bound and best-first node selection, we isolate the effects of the search mechanisms, particularly the identification of equivalent states reached through different scheduling decision sequences. Experiments on 480 benchmark RCPSP instances show that graph-based $A^*$ substantially outperforms the B\&B variants considered. With additional pruning rules, $A^*$ proves optimality for all 480 instances, compared with 445 for best-first B\&B under the same time and memory limits. Depth-first B\&B uses less memory and can find feasible solutions even when the search does not finish within the time limit. Stronger linear programming relaxations reduce the search space but introduce considerable computational overhead. Our findings highlight the importance of duplicate-state detection and clarify the computational trade-offs between the search approaches. These insights can help researchers select an appropriate method and identify opportunities to incorporate effective mechanisms from one approach into the other.
\end{abstract}


\section{Introduction}\label{sec:intro}
The Resource-Constrained Project Scheduling Problem (RCPSP) seeks a schedule that minimizes project completion time while honoring interdependencies among activities and limits on requested resources \citep[see, e.g.,][for a recent overview]{artigues2026fifty}. The problem itself and its extensions model a common challenge in construction \citep[see, e.g.,][]{zh24atd}, manufacturing \citep[see, e.g.,][]{ca12psw}, maintenance \citep[see, e.g.,][]{ch24rcp}, healthcare \citep[see, e.g.,][]{yo26mmr}, software development \citep[see, e.g.,][]{li24sps}, and other project-based settings in which activities that could otherwise run in parallel must compete for scarce resources. Because even moderately sized instances can be difficult to solve optimally, the RCPSP is important both as a practical scheduling problem and as a benchmark for exact optimization methods.

Two mature research traditions approach this problem from different perspectives (details are presented in Section~\ref{sec:lit}). In Operations Research (OR), exact methods commonly explore a tree of restricted subproblems or partial schedules using Branch-and-Bound (B\&B). Lower bounds are used to eliminate subproblems that cannot improve the best complete schedule found so far. In Artificial Intelligence (AI), the problem can instead be represented as a state-space graph and solved using $A^*$, which expands states in best-first order according to an admissible evaluation function and merges equivalent states reached through different decision sequences. At a sufficiently abstract level, $A^*$ can be represented within general formulations of B\&B \citep{nau1984general}. This theoretical relationship, however, does not explain how the concrete algorithms used in OR and AI differ when applied to the same scheduling problem. In practice, B\&B in its \textit{ordinary} variant \citep[see, ordinary B\&B in][]{nau1984general} explores a tree of subsets of feasible solutions implicitly defined by partial schedules, whereas graph-based $A^*$ systematically searches a graph, detecting and merging duplicate states reached through different decision sequences. B\&B may also use the best feasible solution found so far to prune nodes. Existing comparisons are often affected by differences in problem representation, lower-bound strength, node-selection strategy, or algorithm-specific enhancements. Consequently, it remains unclear how much of the observed performance difference is caused by the different search mechanisms themselves.

This study aims to explain the fundamental computational trade-offs between these two widely used exact-search paradigms and to provide a basis for selecting, improving, or combining them. Its purpose is therefore not to determine whether a basic implementation of either method outperforms state-of-the-art RCPSP solvers. Instead, our contributions are threefold. First, we develop a common modeling foundation for the comparison by presenting mixed-integer linear programming (MILP), state-action, and timed Petri-net representations of the RCPSP and establishing the relationships between their schedules, states, and search graphs. Second, we provide a mechanistic and formal comparison of B\&B and \(A^*\) under identical lower-bound information and best-first node selection. The analysis identifies how duplicate states, pruning, repeated heuristic evaluations, and memory requirements affect the two searches. Finally, we conduct a controlled empirical study on benchmark RCPSP instances, comparing computational time, generated and expanded nodes, lower-bound evaluations, and memory use. The experiments examine how instance characteristics and the search mechanisms influence the relative performance of the two methods.

The remainder of the paper is organized as follows. Section~\ref{sec:lit} reviews the relevant literature and identifies the gap addressed by this study. Section~\ref{sec:problem} defines the RCPSP and introduces a running example. Section~\ref{sec:representations} presents the MILP, state-action, and timed Petri-net representations and establishes their connections. Section~\ref{sec:approaches} describes the B\&B and $A^*$ algorithms and compares their search mechanisms.  Section~\ref{sec:empirical} reports the empirical study. Finally, Section~\ref{sec:summary} concludes the paper and discusses directions for future research.

\section{Related Literature} \label{sec:lit}
This section reviews three related research streams: B\&B methods for the RCPSP (Section~\ref{subsec:lit1}), heuristic state-space search for scheduling (Section~\ref{subsec:lit2}), and theoretical studies on the relationship between B\&B and $A^*$ (Section~\ref{subsec:lit3}). Together, these streams position our comparison of optimization-based B\&B and graph-based $A^*$ (Section~\ref{subsec:lit4}).

\subsection{Branch-and-Bound Methods for the RCPSP}\label{subsec:lit1}
Exact RCPSP solution approaches in OR have mainly focused on MILP and specialized B\&B algorithms. Early exact approaches were compared by \citet{patterson1984comparison}. Subsequent B\&B algorithms introduced delay alternatives, schedule schemes, memory-efficient enumeration, and constraint propagation \citep{demeulemeester1992branch,brucker1998branch,sprecher2000scheduling,dorndorf2000branch}. Their performance depends on the branching scheme, lower bounds, dominance rules, node-selection strategy, and the early discovery of good feasible schedules \citep{brucker2003lower,coelho2018exact,guo2023prediction}. In parallel, MILP formulations range from classical time-indexed models \citep{pritsker1969multiproject} to event-based formulations \citep{kone2011event}. Surveys document the continued development of these approaches \citep{brucker1999resource,artigues2026fifty}.

In this OR interpretation, B\&B is an optimization framework. Each node in the search tree represents a subset of feasible schedules, encoded either by a partial schedule or by a restricted optimization subproblem. Branching decisions partition this subset into smaller subsets. A valid lower bound, often but not necessarily obtained from a relaxation, bounds the best schedule attainable within a node, while the best complete schedule found so far provides an upper bound. A node can be discarded if it is infeasible or if its lower bound is not smaller than the upper bound. Node selection is flexible and may follow depth-first, best-first, breadth-first, or hybrid strategies. Equivalent nodes reached through different decision sequences are not systematically merged, although problem-specific dominance rules may remove some of this redundancy.

\subsection{State-Space Search for Scheduling}\label{subsec:lit2}
In $A^*$, nodes are selected according to $\hat f=g+h$, where $g$ is the cost accumulated so far and $h$ is an admissible lower bound on the remaining cost \citep{hart1968formal}. As originally defined, $A^*$ performs graph search: its \textsc{Open} and \textsc{Closed} structures detect when a state is regenerated along a different path, and only the cheapest path to each state is retained. 
This mechanism’s effect depends on the state representation.
Ordinary B\&B, in contrast, may retain the corresponding decision sequences as distinct nodes unless a dominance rule detects their equivalence.

Timed Petri nets have long been used to model scheduling problems, including job-shop systems
\citep{hillion1989petri}, manufacturing systems \citep{proth1998scheduling,lee1994scheduling,mejia2005petri}, and data-driven project planning \citep{cohen2024planning}. Scheduling through reachability-graph search, particularly
using $A^*$, is a well-established approach in this literature \citep{lee1994scheduling,huang2023survey}.

\citet{bell1990solving} developed an exact $A^*$-based algorithm for the RCPSP. Its nodes are subproblems obtained by adding precedence constraints to resolve resource conflicts, so equivalent project states reached through different repair sequences are never identified. Accordingly, \citet{nazareth1999breadth} later interpreted the method as a best-first B\&B procedure formulated as an $A^*$ tree search. More recently, \citet{lublin2026petri} formulated the RCPSP as an optimal search over a Petri-net reachability graph and showed that graph-based $A^*$ can compete with general-purpose MILP solvers. That comparison demonstrates the potential of $A^*$, but it does not isolate the algorithmic differences between $A^*$ and B\&B under identical lower-bound information and search conditions.

\subsection{Relationship between B\&B and \texorpdfstring{$A^*$}{A*}}\label{subsec:lit3}
Earlier studies related B\&B to dynamic programming and state-space representations \citep{morin1976branch,ibaraki1978branch}, while
\citet{kumar1983general} proposed a broad B\&B formulation for AND/OR-tree search. \citet{nau1984general} introduced \emph{General Branch and Bound} (GBB), defined through abstract representation, goal, selection, splitting, and pruning functions. These functions are deliberately general enough for both ordinary B\&B and $A^*$ to be represented as special cases. We argue that this common abstraction should not be interpreted as showing that graph-based $A^*$ and ordinary optimization B\&B are the same practical algorithm. In the terminology of \citet{nau1984general}, ordinary B\&B selects one active subproblem at a time and prunes using a lower bound and an upper bound which is given by the best solution found so far. $A^*$ is embedded in GBB through a different choice of representations, selection rules, and pruning functions. Similarly, \citet{labat2003branch} concluded that heuristic search and B\&B apply the same principles and differ mainly at the interpretation level within an abstract view of B\&B. Thus, the AI literature compares the methods within a highly general search abstraction, rather than comparing the concrete optimization-based B\&B used in RCPSP research with graph-based $A^*$.

\subsection{Research Gap}\label{subsec:lit4}
Taken together, the existing literature does not provide a controlled, practical comparison of the two introduced search paradigms for the RCPSP. An admissible $A^*$ heuristic and a B\&B lower bound serve the same mathematical purpose and may be computed from the same relaxation. Therefore, differences in bound strength should be controlled when comparing the methods. To the best of our knowledge, no study compares ordinary OR-style B\&B and graph-based $A^*$ on identical RCPSP instances while holding the lower-bound function and node-selection conditions fixed. This paper addresses that gap and isolates the consequences of tree-based enumeration and upper-bound pruning in B\&B versus systematic duplicate-state detection in graph-based $A^*$.

\section{The Resource-Constrained Project Scheduling Problem} \label{sec:problem}
A project consists of activities that must be executed subject to precedence relations and limited resources. In the RCPSP, the goal is to choose activity start times that minimize the makespan while satisfying all precedence and resource constraints. This section gives the formal problem definition (Section~\ref{sec:problem_formal}) and introduces the running example (Section~\ref{sec:running_example}) used throughout the paper.

\subsection{Formal Definition}
\label{sec:problem_formal}
Let $\V$ represent the set of non-preemptive activities, with $d_i \in \mathbb{N}_0=\{0,1,2,\ldots\}$ denoting the duration of activity $i \in \V$. Let $\R$ be the set of renewable resources, with capacity $R_k$ for each $k\in\R$. Additionally, let $\E\subset \V\times\V$ be the set of all direct precedence relations, meaning that $(i,j)\in \E$ if the completion of activity $i$ is a direct prerequisite for starting activity $j$. For each activity $i\in\V$, we define $\Pred(i)=\{j\in\V:(j,i)\in\E\}$, $\Suc(i)=\{j\in\V:(i,j)\in\E\}$. The elements of $\Pred(i)$ and $\Suc(i)$ are the immediate predecessors and immediate successors of $i$, respectively.

Without loss of generality, assume that the unique starting point (no predecessors) of the project is at activity $0$ and the unique ending point (no successors) of the project is activity $n$, with $0,n\in \V$ and $d_{0}=d_{n}=0$. We refer to activities $0$ and $n$ as dummy activities. For any other activity $i$, we assume $d_i > 0$. The resource demand of activity $i \in \V$ for resource $k\in\R$ is denoted by $r_{ik}$, where $r_{0k}=r_{nk}=0$ for all $k\in\R$, and the start time of activity $i\in \V$ is denoted by $S_i$. We assume that $(\V,\E)$ is acyclic and that
$0\le r_{ik}\le R_k$ for all $i\in\V$ and $k\in\R$.
For each time $t \in \mathbb{N}_0$, let $\mathcal{A}(t) \subseteq \mathcal{V}$ denote the set of activities in process at time $t$, so that $i \in \mathcal{A}(t)$ if $S_i \leq t < S_i + d_i$. A schedule $S = \{S_i\}_{i\in \V}$ is an assignment of start times to all activities. A schedule is feasible if it satisfies all precedence constraints, and at each point in time the resource demand of the processed activities does not exceed the amount of available resources. A feasible schedule is optimal if its makespan, given by $S_n$, is minimal among all feasible schedules. 

Equivalently, the RCPSP can be stated as the following optimization problem:
\begin{subequations}
\begingroup 
\addtolength{\jot}{-.3em}
\begin{align}
\min_{S}\quad & S_n \\
\text{s.t.}\quad 
& S_i + d_i \le S_j, 
&& \forall (i,j)\in\E,\\
& \sum_{i\in\mathcal{A}(t)} r_{ik} \le R_k,
&& \forall k\in\R,\ t\in\mathbb{N}_0,\\
& S_i \in \mathbb{N}_0,
&& \forall i\in\V .
\end{align}
\endgroup
\end{subequations}

A feasible schedule $S$ is called \emph{active} if there is no feasible schedule $S'\neq S$ that starts every activity no later than in $S$ and starts at least one activity strictly earlier; that is, there does not exist a feasible schedule $S'$ and an activity $j\in\V$ such that $S'_i\leq S_i,\; \forall i\in\V$, and $S'_j<S_j$. Equivalently, $S$ cannot be improved by moving one or more activities earlier without delaying any other activity or violating a precedence or resource constraint. Since at least one optimal schedule is active~\citep{kolisch1996serial},
focusing on earliest-start schedule generation 
preserves optimality. In the search formulations below, once an activity is selected, it is started at the earliest feasible time relative to the current state, subject to precedence and resource feasibility. Hence, the resulting state space contains at least one path corresponding to an optimal active schedule.

We use the term \emph{partial schedule} to refer to a feasible schedule restricted to a subset of already scheduled activities. Formally, a partial schedule is a pair $(\mathcal{J},S_{\mathcal{J}})$, where $\mathcal{J}\subseteq \mathcal{V}$ is predecessor-closed, i.e., if $j\in\mathcal{J}$ and $(i,j)\in\mathcal{E}$, then $i\in\mathcal{J}$, and $S_{\mathcal{J}}=\{S_j\}_{j\in\mathcal{J}}$ assigns start times to all activities in $\mathcal{J}$. The assignment must satisfy all precedence and resource constraints induced by the activities in $\mathcal{J}$. Each partial schedule induces a timed system state that records the current time, which is the start time of the latest scheduled activity, the completed and ongoing activities, and the residual durations of the ongoing activities. This state contains all information required to determine the feasible future scheduling decisions. Consequently, any two partial schedules that induce the same timed state define the same remaining scheduling problem. Our B\&B formulation searches the tree of partial schedules, whereas graph-based $A^*$ searches the corresponding state graph obtained by merging partial schedules that induce the same timed state. Thus, both methods represent the same underlying scheduling decisions but differ in whether equivalent scheduling histories are retained as separate nodes or merged into a single state.

The RCPSP is NP-hard in the strong sense~\citep{blazewicz1983scheduling}, which necessitates the use of effective bounding and pruning strategies within exact methods.

\subsection{Running Example}
\label{sec:running_example}
Throughout this paper, we use a small RCPSP instance to illustrate the concepts and methods. The instance contains the dummy start and end activities $0$ and $n$, respectively, and five non-dummy activities $a,b,c,d,e$. It is illustrated in Figure~\ref{fig:running_example_aon}.

\begin{figure}[t]
\centering
\begin{tikzpicture}[
    aonnode/.style={
        rectangle, draw=black, thick, rounded corners=2pt,
        minimum width=17mm, minimum height=9mm,
        inner sep=2pt, align=center, font=\footnotesize
    },
    dummy/.style={
        rectangle, draw=black, thick, rounded corners=2pt,
        minimum width=17mm, minimum height=9mm,
        inner sep=2pt, align=center, font=\footnotesize,
        fill=gray!10
    },
    arc/.style={-{Stealth[length=5pt,width=4pt]}, thick},
    ann/.style={font=\scriptsize, align=center}
]

\node[dummy] (v0) at (0,0)
{$0$\\{\scriptsize $d_0=0,\ r_0=0$}};

\node[aonnode] (a) at (3.0, 1.6)
{$a$\\{\scriptsize $d_a=3,\ r_a=1$}};
\node[aonnode] (b) at (3.0, 0.0)
{$b$\\{\scriptsize $d_b=2,\ r_b=2$}};
\node[aonnode] (c) at (3.0,-1.6)
{$c$\\{\scriptsize $d_c=3,\ r_c=1$}};

\node[aonnode] (d) at (6.2, 0.8)
{$d$\\{\scriptsize $d_d=1,\ r_d=1$}};
\node[aonnode] (e) at (6.2,-0.8)
{$e$\\{\scriptsize $d_e=1,\ r_e=1$}};

\node[dummy] (vn) at (9.2,0)
{$n$\\{\scriptsize $d_n=0,\ r_n=0$}};

\draw[arc] (v0) -- (a);
\draw[arc] (v0) -- (b);
\draw[arc] (v0) -- (c);

\draw[arc] (a) -- (d);
\draw[arc] (b) -- (d);

\draw[arc] (b) -- (e);
\draw[arc] (c) -- (e);

\draw[arc] (d) -- (vn);
\draw[arc] (e) -- (vn);

\end{tikzpicture}
\caption{Illustration of the running RCPSP example.}
\vspace{-.3cm}\caption*{\scriptsize \textit{Note:} Each node shows the activity identifier, duration $d$, and demand $r$ for the single renewable resource with capacity $R_1=2$.}
\label{fig:running_example_aon}
\end{figure}

The project comprises the seven activities $\V=\{0,a,b,c,d,e,n\}$. Activities $0$ and $n$ have zero duration. The durations of the non-dummy activities are $d_a=3$, $d_b=2$, $d_c=3$, $d_d=1$, and $d_e=1$. The precedence relations are completely captured by $\Suc(0)=\{a,b,c\}$, $\Pred(d)=\{a,b\}$, $\Pred(e)=\{b,c\}$, and $\Pred(n)=\{d,e\}$. Equivalently, using informal arrow notation, they can be written as $0\to\{a,b,c\}$, $\{a,b\}\to d$, $\{b,c\}\to e$, and $\{d,e\}\to n$. A single renewable resource with capacity $R_1=2$ is available. The resource demands are $r_a=1$, $r_b=2$, $r_c=1$, $r_d=1$, and $r_e=1$, while the dummy activities $0$ and $n$ require no resources. Since the example contains a single resource, we omit the resource index; for example, $r_{a1}$ is written as $r_a$. The example is small enough to be solved by hand, yet rich enough to exhibit key differences between B\&B and $A^*$ search. Concurrent activities create multiple equivalent execution orderings, while the resource constraint prevents full parallelism and forces nontrivial scheduling decisions.

\section{Problem Representations}
\label{sec:representations}
The RCPSP can be modeled from several computational perspectives. A MILP formulation (Section~\ref{sec:milp}) represents complete schedules through decision variables and supports ordinary OR-style B\&B. A state-action formulation (Section~\ref{sec:dp}) represents scheduling as transitions between timed system states and supports graph-search methods such as $A^*$. A Petri-net formulation (Section~\ref{sec:petri_nets}) provides an explicit token-based model of precedence, timing, and renewable-resource availability, whose timed reachability graph is equivalent to the state graph. This section introduces the three formulations and formalizes their relationships (Section~\ref{sec:equivalences}).

\subsection{Mixed-Integer Linear Programming Formulation}
\label{sec:milp}
The RCPSP can be formulated as a MILP by introducing binary decision variables that encode the scheduling decisions. Several MILP formulations exist in the literature; here we present a standard time-indexed formulation.

Let $H$ denote an upper bound on the project duration, usually obtained by scheduling all activities in sequence so $H=\sum_{i\in \V} d_i$, and let $\T=\{0,1,\ldots,H\}$ denote the corresponding discrete time horizon.
Define binary variables
\begin{equation}
y_{it} =
\begin{cases}
1, & \text{if activity } i \text{ starts at time } t,\\[-.3cm]
0, & \text{otherwise,}
\end{cases}
\end{equation} for all $i \in \V$ and $ t\in \T$.
Then, the MILP formulation is as follows.
\begin{subequations} \label{rcpsp_MILP} \allowdisplaybreaks
\begingroup 
\addtolength{\jot}{-.3em}
\begin{align}
\min\ &S_{n}\label{rcpsp_obj}
\\
\text{s.t.}\quad
& \sum_{t=0}^{H-d_i} y_{it} = 1
&& \forall i \in \V\label{rcpsp_once}\\
& S_i = \sum_{t=0}^{H-d_i} t\, y_{it}
&& \forall i \in \V\label{rcpsp_start}\\
& S_i + d_i \le S_j 
&& \forall (i,j) \in \E\label{rcpsp_prec}\\
& \sum_{i \in \V} r_{ik}
\sum_{\tau=\max\{0,t-d_i+1\}}^{\min\{t,H-d_i\}} y_{i\tau}
\le R_k
&& \forall t \in \T,\  k \in \R\label{rcpsp_res}\\
& S_i \in \mathbb{N}_0
&& \forall i \in \V\label{rcpsp_int}\\
& y_{it}\in\{0,1\}
&& \forall i \in \V,\ t=0,\dots,H-d_i\label{rcpsp_binary}
\end{align}
\endgroup
\end{subequations}

The objective~\eqref{rcpsp_obj} minimizes the makespan, represented by the start time of the dummy activity $n$, $S_n$. Constraints~\eqref{rcpsp_once} ensure that each activity starts exactly once, while constraints~\eqref{rcpsp_start} define the corresponding start times $S_i$.
Constraints~\eqref{rcpsp_prec} enforce the precedence relations. For a fixed activity $i$ and time period $t$, activity $i$ is ongoing during time period $t$ if and only if it started at some time $\tau$ with $t-d_i+1 \le \tau \le t$. Since start variables are defined only for feasible start times $\tau \in \{0,\dots,H-d_i\}$, the inner summation in~\eqref{rcpsp_res} ranges over the intersection of these two intervals, namely $\tau=\max\{0,t-d_i+1\},\dots,\min\{t,H-d_i\}$.  Finally, constraints~\eqref{rcpsp_int} and~\eqref{rcpsp_binary} impose integrality on the start-time and binary decision variables.

\subsection{State-Action  Representation}
\label{sec:dp}
An alternative representation views the RCPSP as a shortest-path problem in a state space, where each state encodes the current execution status of the project and transitions (actions) correspond to scheduling decisions. More specifically, a state $s$ captures the project status at a particular point in time. It records:

\begin{itemize}
     \item The current system time $g(s) \in \mathbb{N}_0$, i.e., the time represented by state $s$ at which the system is observed and the next scheduling decision is considered.
    \item The set of completed activities $\V^c(s) \subseteq \V$, i.e., activities that have finished no later than $g(s)$.
    \item The set of ongoing activities $\V^o(s) \subseteq \V$ that have started but have not yet completed at $g(s)$, together with their residual durations $\tilde{d}_i(s) \in \mathbb{N}_0$ relative to $g(s)$ for all $i \in \V^o(s)$, representing the remaining processing time at the current system time $g(s)$. 
\end{itemize}
Although $\V^o(s)$ is a set, in the state tuple we represent its elements as an ordered tuple, with activities listed in a fixed order and their corresponding residual durations in $\tilde d(s)$ listed in the same order.
Activities in $\V^u(s)=\V\setminus\bigl(\V^c(s)\cup\V^o(s)\bigr)$ are called not-yet-started activities.

We define an initial state $s_0$ with $\V^c(s_0) = \emptyset$, $\V^o(s_0) = \emptyset$, $g(s_0) = 0$, and full resource availability. Furthermore, we define a set of goal states $\Sset_f$ as the set of all states $s$ with  $\V^c(s) = \V$. The available resource profile at time $g(s)$ can be deduced from resources consumed by ongoing activities and their capacities. A state $s$ is called feasible if it satisfies the resource constraints $\sum_{i\in \V^o(s)} r_{ik}\leq R_k$ for all $k\in\R$, and the precedence constraints, i.e., if $j\in \V^o(s)$ and $(i,j)\in\E$, then $i\in \V^c(s)$.

At each state $s$, we can select an activity to be scheduled next, which leads to a new state $s'$. The set of selectable activities at state $s$, denoted by $\V^a(s)$, consists of the not-yet-started activities whose immediate predecessors have either completed or are currently ongoing: $\V^{a}(s)= \left\{ i\in\V^u(s): \Pred(i)\subseteq \V^c(s)\cup\V^o(s) \right\}$. Thus, $\V^a(s)\subseteq \V^u(s)$. If some predecessors of the selected activity are still ongoing, or if insufficient resources are currently available, the system time advances to the earliest time at which the selected activity can start. This earliest-start interpretation is used in the remainder of the paper. Hence, for a given activity $j\in \V^{a}(s)$, we get:
\begin{itemize}
\item Time transition:
\begin{equation}\label{eq:time_transition}
g(s')=\min\left\{t:
\begin{array}{l}
t\geq g(s),\\[1pt]
t\geq g(s)+\tilde{d}_i(s)
\quad \forall i\in \V^o(s)\ \text{with } (i,j)\in\E,\\[1pt]
r_{jk} \leq R_k-\displaystyle\sum_{i\in \V^o(s):\,g(s)+\tilde{d}_i(s)> t} r_{ik},
\quad \forall k\in\R
\end{array}
\right\},
\end{equation}
with $\Delta(s,s')=g(s')-g(s)$ denoting the elapsed time in the transition. 

\item Ongoing activities:
\begin{equation}\label{eq:ongoing_transition}
\V^o(s')=
\{i\in \V^o(s): \Delta(s,s')<\tilde{d}_i(s)\}
\cup
\begin{cases}
\{j\}, & \text{if } d_j>0,\\[-.2cm]
\emptyset, & \text{if } d_j=0.
\end{cases}
\end{equation}
Thus, an activity is ongoing in $s'$ 
when it has positive residual duration after the transition. Residual durations of previously ongoing activities are reduced by the elapsed time, while the newly started activity $j$ receives its full duration if it is positive. Hence, for all $i\in \V^o(s')\setminus\{j\}$, we have $\tilde{d}_i(s')=\tilde{d}_i(s)-\Delta(s,s')$, and, if $d_j>0$, $\tilde{d}_j(s')=d_j$.

\item Completed activities:
\begin{equation}\label{eq:completed_transition}
\V^c(s')=
\V^c(s)
\cup
\{i\in \V^o(s): \Delta(s,s')\geq \tilde{d}_i(s)\}
\cup
\begin{cases}
\{j\}, & \text{if } d_j=0,\\[-.2cm]
\emptyset, & \text{if } d_j>0.
\end{cases}
\end{equation}

\end{itemize}

The reachable state space can be represented as a directed graph $G=(\Sset,\mathbb{E})$, which we call the state graph, where $\Sset$ is the set of states reachable from $s_0$ by feasible scheduling decisions, and $\mathbb{E}$ is the set of directed edges induced by these decisions. We write $ s \xrightarrow{\,j\,} s'$ if choosing activity $j$ in state $s$ leads to the successor state $s'$ under the above transition rules. Finding an optimal schedule is equivalent to finding a path from the initial state to a goal state that minimizes the system time of the goal state. Thus, the cost of a path is measured by the system time at its final state. This equivalence follows because an optimal schedule can be represented by ordering its activities according to their start times and applying the earliest-start transition rule along that order.

To illustrate this notation, consider the running example. We start from state $s_0$. Since $\V^a(s_0)=\{0\}$, scheduling dummy activity $0$ leads to $s_1=\bigl(g(s_1),\V^c(s_1),\V^o(s_1),\tilde d(s_1)\bigr)=(0,\{0\},(),())$. At this state, the set of selectable activities is $\V^a(s_1)=\{a,b,c\}$. If activity $a$ is selected next, the search moves to $s_2=(0,\{0\},(a),(3))$, and the selectable activities become $\V^a(s_2)=\{b,c\}$. If activity $b$ is then selected, the transition advances time to state $s_5$. Activity $b$ requires two resource units, one of which is currently used by activity $a$, and therefore cannot start until $a$ releases that unit. Consequently, $g(s_5)=3$ and $s_5=(3,\{0,a\},(b),(2))$. In contrast, if activity $c$ is selected after activity $a$, then $a$ and $c$ can be processed in parallel because sufficient resource capacity is available. The resulting state is $s_6=(0,\{0\},(a,c),(3,3))$. 
The same state $s_6$ is reached by scheduling $c$ before $a$, illustrating how different scheduling decision sequences can lead to the same state.
Figure~\ref{fig:rcpsp-state-space} illustrates the subgraph of $G$ induced by these scheduling decisions.
\begin{figure}[htbp]
    \centering
\begin{tikzpicture}[
    state/.style={
        rectangle, rounded corners=2pt,
        draw=black, thick, fill=white,
        minimum width=30mm,
        minimum height=11mm,
        align=center,
        font=\footnotesize,
        inner sep=2pt
    },
    arc/.style={-{Stealth[length=5pt,width=4pt]}, thick},
    tlbl/.style={font=\footnotesize, fill=white, inner sep=1.5pt},
]


\node[state] (S0) at (0,0)
{$s_0$\\
$(0,\varnothing,(),())$};

\node[state] (S1) at (4.0,0)
{$s_1$\\
$(0,\{0\},(),())$};

\node[state] (S2) at (8.3,1.7)  {$s_2$\\ $(0,\{0\},(a),(3))$};
\node[state] (S3) at (8.3,0)    {$s_3$\\ $(0,\{0\},(b),(2))$};
\node[state] (S4) at (8.3,-1.7) {$s_4$\\ $(0,\{0\},(c),(3))$};
\node[state] (S5) at (12.8,1.7) {$s_5$\\ $(3,\{0,a\},(b),(2))$};
\node[state] (S6) at (12.8,-1.7){$s_6$\\ $(0,\{0\},(a,c),(3,3))$};

\draw[arc] (S0) -- node[tlbl, above] {$0$} (S1);
\draw[arc] (S1) -- node[tlbl, above] {$a$} (S2);
\draw[arc] (S1) -- node[tlbl, above, pos=0.6] {$b$} (S3);
\draw[arc] (S1) -- node[tlbl, below] {$c$} (S4);
\draw[arc] (S2) -- node[tlbl, above] {$b$} (S5);
\draw[arc] (S2.south east) -- node[tlbl, right] {$c$} (S6.north west);
\draw[arc] (S4) -- node[tlbl, below] {$a$} (S6);

\end{tikzpicture}
    \caption{Illustrative subgraph of the RCPSP state graph for the first scheduling decisions.}
    \label{fig:rcpsp-state-space}
    \vspace{-.3cm}\caption*{\scriptsize \textit{Note:} Each node is a state represented as $(g(s),\V^c(s),\V^o(s),\tilde{d}(s))$, where $g(s)$ is the current time, $\V^c$ is the set of completed activities, $\V^o$ is the ordered tuple of ongoing activities, and $\tilde{d}$ is the corresponding ordered tuple of residual processing times.}
\end{figure}

\subsection{Petri Net Representation}
\label{sec:petri_nets}
A Petri net can be viewed as a directed bipartite graph with two types of nodes: places and transitions. 
Petri nets provide a formal mathematical framework for modeling concurrent systems, capturing synchronization, precedence relations, and shared resources \citep{murata1989petri}. We use them to model concurrent project activities and resource availability, and to establish a formal correspondence between their timed reachability graph and the state graph used by $A^*$.
Places represent conditions, intermediate project states, or available resources. Transitions represent activities or events. Directed arcs connect only places to transitions or transitions to places. Therefore, there are no direct arcs between two places or between two transitions. A token is a marker residing in a place, and the distribution of tokens over the places represents the current state of the system. The movement of tokens through the net naturally represents synchronization and parallelism among activities. Petri nets can also represent XOR choices between alternative courses of action. The RCPSP considered here, however, contains only AND-type precedence relations. An activity may start only after all of its predecessors have completed. Accordingly, we do not use XOR constructs in the Petri-net representation. In our comparison, the Petri-net representation provides an established modeling framework for state-based search, rather than a requirement for applying $A^*$.

Petri nets support correctness checks such as whether all activities of a project are executable, the project can always progress to completion, and no residual work remains in the final state. They thus provide verification capabilities that go well beyond what alternative representations are capable of \citep{murata1989petri,van2022process}. Moreover, Petri nets can be used to model and visualize resource constraints.

\subsubsection{Timed Transition Petri Nets with Resources}
\label{sec:TTPNR_def}
We model the RCPSP as a so-called Timed Transition Petri Net with Resources ($TTPNR$), i.e., a tuple $TTPNR = (P, T, F, W, M_0, L, \tau, \mathcal{R})$, where, in addition to the definitions in Section~\ref{sec:problem_formal}, we define: 

\begin{itemize}
    \item $P=P_C\cup P_R$ is the finite set of \textit{places}, depicted as circles, where $P_C\cap P_R=\emptyset$. The set $P_C$ contains the control-flow places $\{p_1, p_2, \ldots, p_m\}$, which represent the progress of the project. The set $P_R=\{p_{r_k}:k\in\R\}$ contains one resource place for each renewable resource. Tokens in resource places represent available resource units.  
    \item $T=\{t_i:i\in\V\}$ is the finite set of \textit{transitions}, depicted as rectangles, with $P\cap T=\emptyset$. Transition $t_i$ represents activity $i$. In particular, the dummy start and end activities are represented by zero-duration transitions $t_0$ and $t_n$, respectively.
    \item $F \subseteq (P \times T) \cup (T \times P)$ is the set of directed arcs, also called the \textit{flow relation}. We distinguish between control-flow arcs and resource arcs. 
    Control-flow arcs connect control-flow places and transitions and encode the precedence/progress structure of the project. 
    Resource arcs connect resource places and transitions and model the use of renewable resource units. For the RCPSP, the control-flow structure consists of a source place $p_{\mathrm{start}}$ connected to $t_0$, a sink place $p_{\mathrm{end}}$ connected from $t_n$, and one place $p_{ij}$ connecting $t_i$ to $t_j$ for each $(i,j)\in\E$. These are the only control-flow arcs, and initially only $p_{\mathrm{start}}$ contains a token. In the running example (Figures \ref{fig:toyTTPNR} and \ref{fig_A_star_on_RG}), these places are denoted by $p_1,\ldots,p_{11}$.

    \item $W : F \rightarrow \mathbb{Z}^+$ is a weight function that assigns a positive integer weight to each arc. For an input arc $(p,t)$, the value $W(p,t)$ is the number of tokens that must be available in place $p$ for transition $t$ to start. For an output arc $(t,p)$, the value $W(t,p)$ is the number of tokens deposited in place $p$ after transition $t$ completes. In our RCPSP models, control-flow arcs have weight one. Resource arcs have weight equal to the number of resource units required by the corresponding activity, and the same number of resource tokens is returned to the resource place when the activity completes.
    \item A marking $M:P\rightarrow\mathbb{N}_0$ specifies the number of
    tokens in each place and thereby encodes the current system state.
    Specifically, $M_0:P\rightarrow\mathbb{N}_0$ is the initial marking.
    For control-flow places, $M_0$ encodes the initial project state. For each resource place $p_{r_k}\in P_R$, $M_0(p_{r_k})=R_k$ specifies the initial number of available units of resource $k$. In the timed state
    representation introduced in Section~\ref{sec:RG_TTPNR}, this ordinary
    marking is refined by associating a residual availability time with each token, defined as the number of time units, measured from the current system time, until the token becomes available. The underlying
    distribution of tokens over places remains the marking defined here.
    \item $L:T\rightarrow\V$ is a labeling function that assigns each transition its corresponding project activity. Since transition $t_i$ represents activity $i$, we have $L(t_i)=i$ for all $i\in\V$.
    \item $\tau:T\rightarrow\mathbb{N}_0$ is a timing function that specifies the duration of each transition. Since transition $t_i$ represents activity $i$, we set $\tau(t_i)=d_i$ for all $i\in\V$. Starting, or firing, a transition means consuming the required input tokens and committing to the corresponding activity. When transition $t_i$ fires, the required input tokens are consumed, and its output control-flow tokens and returned resource tokens are associated with their corresponding places with residual availability time $\tau(t_i)$. These tokens can enable subsequent transitions only when their residual availability time reaches zero.
\end{itemize}
When no resources are present ($\mathcal{R} = \emptyset$, $P_R = \emptyset$) and all durations are zero ($\tau(t) = 0$ for all $t$), the $TTPNR$ reduces to a standard Petri Net. When resources are absent but timing is present, it reduces to a Timed Transition Petri Net.

For compactness, we often represent markings by listing only places with positive token counts. Thus, $[p_1]$ denotes the marking $M$ with $M(p_1)=1$ and $M(p)=0$ for all $p\neq p_1$. More generally, $p^q$ denotes $q$ tokens in place $p$, while the exponent is omitted when $q=1$. For example, $[p_1,p_{r_1}^2]$ denotes a marking with one token in the control-flow place $p_1$, two tokens in the resource place $p_{r_1}$, and zero tokens in all other places.

A place $p$ is an input place of transition $t$ if $(p, t) \in F$, and an output place if $(t, p) \in F$. We write $\bullet t = \{p \in P : (p, t) \in F\}$ and $t \bullet = \{p \in P : (t, p) \in F\}$ for the sets of input and output places of $t$, respectively. Similarly, $\bullet p$ and $p \bullet$ denote the input and output transitions of place $p$.

A transition $t \in T$ is enabled in marking $M$ if $M(p) \geq W(p, t)$ for all $p \in \bullet t$. This condition requires that every input place, including both control-flow and resource places, contains sufficiently many tokens. 

When $t$ fires, it consumes tokens from its input places and produces tokens for its output places. The corresponding structural marking update  $M'$ is given by:
\begin{equation}
M'(p)
=
M(p)
-
\begin{cases}
W(p,t), & (p,t)\in F,\\[-.2cm]
0, & \text{otherwise,}
\end{cases}
+
\begin{cases}
W(t,p), & (t,p)\in F,\\[-.2cm]
0, & \text{otherwise,}
\end{cases}
\qquad \forall p\in P .
\label{eq:marking}
\end{equation}

Equation~\eqref{eq:marking} describes the structural change in the distribution of tokens over places. In the timed interpretation used in this paper, however, the presence of a token in a place does not necessarily imply that it is immediately available. The availability of produced control-flow tokens and returned resource tokens depends on the duration $\tau(t)$. This timing information is incorporated into the state representation introduced in Section~\ref{sec:RG_TTPNR}.

\begin{figure}[t]
\centering
\begin{tikzpicture}[
    place/.style={circle, draw, minimum size=6mm, inner sep=0pt},
    resource/.style={circle, draw, black!50, thick, minimum size=6mm, inner sep=0pt},
    transition/.style={rectangle, draw, minimum width=6mm, minimum height=6mm},
    arc/.style={->, >=stealth, thick},
    redarc/.style={->, >=stealth, thick, dashed,  black!50},
    scale=0.85, transform shape
]
    
    \node[place, tokens=1, label=below:$p_1$] (p1) at (0,0) {};
    \node[place] (p2) at (2.7,1.2) {$p_2$};
    \node[place] (p3) at (2.7,0) {$p_3$};
    \node[place] (p4) at (2.7,-1.2) {$p_4$};
    \node[place] (p5) at (6,1.2) {$p_5$};
    \node[place] (p6) at (5,0.4) {$p_6$};    
    \node[place] (p7) at (5,-0.4) {$p_7$};   
    \node[place] (p8) at (6,-1.2) {$p_8$};
    \node[place] (p9) at (8,1.2) {$p_9$};
    \node[place] (p10) at (8,-1.2) {$p_{10}$};
    \node[place] (p11) at (11,0) {$p_{11}$};
    
    \node[resource, label=right:\textcolor{black!50}{$p_{r_1}$}] (pr) at (7.0,0) {};
    \fill[black!50] (6.9,0.1) circle (1.5pt);
    \fill[black!50] (7.1,0.1) circle (1.5pt);
    
    
    \node[transition, label={below:$\tau=0$}](t1) at (1.0,0) {$t_0$};
    \node[transition, label={above:$\tau=3$}](t2) at (4,1.2) {$t_a$};
    \node[transition, label={below:$\tau=2$}](t3) at (4,0) {$t_b$};
    \node[transition, label={below:$\tau=3$}](t4) at (4,-1.2) {$t_c$};
    \node[transition, label={above:$\tau=1$}](t5) at (7,1.2) {$t_d$};
    \node[transition, label={below:$\tau=1$}](t6) at (7,-1.2) {$t_e$};
    \node[transition, label={below:$\tau=0$}](t7) at (10,0) {$t_n$};
    
    \draw[arc] (p1) -- (t1);
    \draw[arc] (t1) -- (p2);
    \draw[arc] (t1) -- (p3);
    \draw[arc] (t1) -- (p4);
    \draw[arc] (p2) -- (t2);
    \draw[arc] (t2) -- (p5);
    \draw[arc] (p3) -- (t3);
    \draw[arc] (t3) -- (p6);
    \draw[arc] (t3) -- (p7);
    \draw[arc] (p4) -- (t4);
    \draw[arc] (t4) -- (p8);
    \draw[arc] (p5) -- (t5);
    \draw[arc] (p6) -- (t5);
    \draw[arc] (t5) -- (p9);
    \draw[arc] (p7) -- (t6);
    \draw[arc] (p8) -- (t6);
    \draw[arc] (t6) -- (p10);
    \draw[arc] (p9) -- (t7);
    \draw[arc] (p10) -- (t7);
    \draw[arc] (t7) -- (p11);
    
    \draw[redarc, transform canvas={yshift=1mm}] (pr) to (t2);
    \draw[redarc, transform canvas={yshift=-0mm}] (t2) to (pr);
    
   \draw[redarc, transform canvas={yshift=0.5mm}] (pr) -- node[above, black!50, pos=0.5, font=\scriptsize, inner sep=1pt] {2} (t3);
    \draw[redarc, transform canvas={yshift=-0.5mm}] (t3) -- node[below, black!50, pos=0.5, font=\scriptsize, inner sep=1pt] {2} (pr);
    
    \draw[redarc, transform canvas={yshift=-1mm}] (pr) to (t4);
    \draw[redarc, transform canvas={yshift=0mm}] (t4) to (pr);
    
    \draw[redarc, transform canvas={xshift=-0.5mm}] (pr) -- (t5);
    \draw[redarc, transform canvas={xshift=0.5mm}] (t5) -- (pr);
    
    \draw[redarc, transform canvas={xshift=-0.5mm}] (pr) -- (t6);
    \draw[redarc, transform canvas={xshift=0.5mm}] (t6) -- (pr);
\end{tikzpicture}
\caption{Timed transition Petri net with resources for the running RCPSP example. }
\label{fig:toyTTPNR}
\vspace{-.3cm}\caption*{\scriptsize \textit{Note:} Black solid arcs encode control-flow structure and dashed gray arcs encode consumption and release of the renewable resource. The resource place $p_{r_1}$ initially contains two tokens, corresponding to capacity $R_1=2$. Each transition is labeled by its transition name and duration $\tau$.}
\end{figure}

In Figure~\ref{fig:toyTTPNR}, we instantiate the $TTPNR$ for the running RCPSP example with $P=\{p_1,\ldots,p_{11},p_{r_1}\}$ and $T=\{t_0,t_a,t_b,t_c,t_d,t_e,t_n\}$. The black solid arcs encode the control-flow structure, i.e., the precedence relations between activities, while the dashed gray arcs encode consumption and release of the single renewable resource with capacity $R_1=2$. The notation directly follows the general activity-indexed transition notation $t_i$, $i\in\V$. Thus, $t_0$ and $t_n$ represent the dummy start and end activities, respectively, while $t_a,\ldots,t_e$ represent the five non-dummy activities of the running example. For every activity $i\in\{0,a,b,c,d,e,n\}$, the labeling and timing functions satisfy $L(t_i)=i$  and $\tau(t_i)=d_i$. The control-flow places encode the precedence structure of the running example. For instance, transition $t_d$, corresponding to activity $d$, requires tokens in both $p_5$ and $p_6$, thereby enforcing that activities $a$ and $b$ must complete before $d$ can start. Similarly, transition $t_e$, corresponding to activity $e$, requires tokens from the output places of activities $b$ and $c$. The resource place $p_{r_1}$ initially contains two tokens, represented by the initial marking $M_0=[p_1,p_{r_1}^2]$. Resource arcs connect $p_{r_1}$ to the non-dummy activity transitions $t_a,t_b,t_c,t_d,t_e$. Their weights satisfy $W(p_{r_1},t_b)=W(t_b,p_{r_1})=2$, reflecting that activity $b$ requires two resource units, whereas each of the other non-dummy activities requires one resource unit.

\subsubsection{State Representation and Search Graph Construction}\label{sec:RG_TTPNR}

The dynamic behavior of a $TTPNR$ can be represented by a reachability graph $RG=(\Sset_{RG},\mathbb{E}_{RG})$. The node set $\Sset_{RG}$ consists of reachable timed states. Each node $s\in\Sset_{RG}$ is represented by a current system time $g(s)$ and a timed marking $M_s$. The edge set $\mathbb{E}_{RG}$ contains directed edges induced by transition firings, where each edge corresponds to firing one transition, i.e., to one scheduling decision. A sequence of transition firings is called a trace and represents one possible project realization sequence. If several transitions can be selected, each such transition defines a different successor edge.

The timed marking $M_s$ refines the ordinary marking $M:P\rightarrow\mathbb{N}_0$ defined in Section~\ref{sec:TTPNR_def}. It preserves the distribution of tokens over places while associating a residual availability time with each token. Formally, suppose that place $p$ contains $q$ tokens in the ordinary marking underlying state $s$. Let $\nu_p$ denote the number of distinct residual availability times of these tokens, $\theta_{p,1},\ldots,\theta_{p,\nu_p}$ denote these times, and $q_{p,\ell}$ denote the number of tokens with residual availability time $\theta_{p,\ell}$. The timed marking at $p$ is represented as \[M_s(p)=\left\{q_{p,1}^{\theta_{p,1}}, \ldots, q_{p,\nu_p}^{\theta_{p,\nu_p}} \right\},\qquad \sum_{\ell=1}^{\nu_p}q_{p,\ell}=q.\] Here, $q_{p,\ell}^{\theta_{p,\ell}}$ denotes $q_{p,\ell}$ tokens that become available after $\theta_{p,\ell}$ time units. If $q=0$, then $M_s(p)=\emptyset$. At the initial state $s_0$, every token is immediately available. Thus, if $M_0(p)=q>0$, then $M_{s_0}(p)=\{q^0\}$.

For compactness, we list only places containing at least one token. Tokens having the same residual availability time are grouped, and a superscript equal to zero is omitted. Thus, $p:2^\theta$ denotes two tokens in place $p$ that both become available after $\theta$ time units, whereas $p:1,1^3$ denotes one immediately available token and one token that becomes available after three time units. In particular, $M_{s_0}(p)=\{q^0\}$ is written as $M_{s_0}(p)=\{q\}$. For example, $M_s=[p_1:1,\ p_{r_1}:1,1^3]$ represents one immediately available token in place $p_1$ and two tokens in resource place $p_{r_1}$, one immediately available and one available after three time units.

\begin{figure}[t]
\centering
\resizebox{\textwidth}{!}{\begin{tikzpicture}[
    state/.style={
        circle, draw=black, thick, fill=white,
        minimum size=10mm, inner sep=0pt, font=\footnotesize
    },
    gstate/.style={
        circle, draw=black, thick, fill=gray!30,
        minimum size=10mm, inner sep=0pt, font=\footnotesize
    },
    goalstate/.style={
        circle, draw=black, thick, fill=gray!30,
        double, double distance=1.2pt,
        minimum size=10mm, inner sep=0pt, font=\footnotesize
    },
    arc/.style={-{Stealth[length=5pt,width=4pt]}, thick},
    dasharc/.style={-{Stealth[length=5pt,width=4pt]}, thick, dashed},
    tlbl/.style={font=\footnotesize, fill=white, inner sep=1.5pt},
    info/.style={font=\scriptsize, align=left},
    fval/.style={font=\scriptsize},
]

\node[gstate]    (S0) at (0, 0)      {$s_0$};
\node[gstate]    (S1) at (3.3, 0)    {$s_1$};
\node[gstate]    (S2) at (7.2, 0)    {$s_2$};
\node[gstate]    (S6) at (13.0, 0)   {$s_6$};

\node[gstate]    (S4) at (4.5, -1.8) {$s_4$};

\node[gstate]    (S3) at (2.0, -5.2) {$s_3$};
\node[state]     (S5) at (7.8, -4.2) {$s_5$};
\node[gstate] (S7) at (12.0, -2.8){$s_{10}$};

\node[state] (S8) at (0.5, -8.0) {$s_8$};
\node[state] (S9) at (3.5, -8.0) {$s_7$};
\node[state] (S10) at (6.2,-8.0){$s_9$};
\node[gstate] (S12) at (10.2,-8.0){$s_{12}$};
\node[gstate]    (S11) at (14.0,-4.0){$s_{11}$};

\node[gstate]    (S13) at (14.0,-6.0){$s_{13}$};

\node[goalstate] (S14) at (14.0,-8.0){$s_{14}$};


\node[fval, anchor=south] at ([xshift=-6pt,yshift=6pt]S0.north)
  {$(0,4,4)$};
\node[info, anchor=south west] at ([xshift=-3pt,yshift=11pt]S0.north east)
  {$[p_1\!:\!1,\; p_{r_1}\!:\!2]$};

\node[fval, anchor=south] at ([yshift=6pt]S1.north) {$(0,4,4)$};
\node[info, anchor=south west] at ([xshift=6pt,yshift=6pt]S1.north east)
  {$[p_2\!:\!1,\; p_3\!:\!1,$\\$ p_4\!:\!1,\;p_{r_1}\!:\!2]$};

\node[fval, anchor=south] at ([yshift=6pt]S2.north) {$(0,4,4)$};
\node[info, anchor=south west] at ([xshift=6pt,yshift=6pt]S2.north east)
  {$[p_3\!:\!1,\; p_4\!:\!1,\;p_5\!:\!1^3,$\\$p_{r_1}\!:\!1,\;1^3]$};

\node[fval, anchor=south] at ([yshift=6pt]S6.north) {$(0,4,4)$};
\node[info, anchor=south west] at ([xshift=6pt,yshift=6pt]S6.north east)
  {$[p_3\!:\!1,\; p_5\!:\!1^3,$\\$p_8\!:\!1^3,\;p_{r_1}\!:\!2^3]$};

\node[fval, anchor=south] at ([xshift=10pt,yshift=2pt]S4.north) {$(0,4,4)$};
\node[info, anchor=north west] at ([xshift=-8pt,yshift=-6pt]S4.east)
  {$[p_2\!:\!1,\; p_3\!:\!1,$\\$p_8\!:\!1^3,\;p_{r_1}\!:\!1,\;1^3]$};

\node[fval, anchor=south east] at ([xshift=5pt,yshift=3pt]S3.north west) {$(0,4,4)$};
\node[info, anchor=east] at ([xshift=102pt,yshift=28pt]S3.west)
  {$[p_2\!:\!1,\; p_4\!:\!1,\;p_6\!:\!1^2,$\\$p_7\!:\!1^2,\; p_{r_1}\!:\!2^2]$};

\node[fval, anchor=east] at ([yshift=7pt]S5.north) {$(3,4,7)$};
\node[info, anchor=north] at ([xshift=5pt, yshift=-2pt]S5.south)
  {$[p_4\!:\!1,\; p_5\!:\!1,\;p_6\!:\!1^2,$\\$p_7\!:\!1^2,\;p_{r_1}\!:\!2^2]$};

\node[fval, anchor=south west] at ([xshift=6pt,yshift=3pt]S7.north east) {$(3,3,6)$};
\node[info, anchor=east] at ([xshift=10pt,yshift=25pt]S7.west)
  {$[p_5\!:\!1,\; p_6\!:\!1^2,$\\$p_7\!:\!1^2,\; p_8\!:\!1,\; p_{r_1}\!:\!2^2]$};

\node[fval, anchor=south east] at ([xshift=15pt,yshift=3pt]S8.north west) {$(2,4,6)$};
\node[info, anchor=north] at ([yshift=-4pt]S8.south)
  {$[p_2\!:\!1,\; p_6\!:\!1,\;p_7\!:\!1,$\\$p_8\!:\!1^3,\; p_{r_1}\!:\!1,\;1^3]$};

\node[fval, anchor=south] at ([xshift=20pt,yshift=3pt]S9.north west) {$(2,4,6)$};
\node[info, anchor=north] at ([yshift=-4pt]S9.south)
  {$[p_4\!:\!1,\; p_5\!:\!1^3,\; p_6\!:\!1,$\\$p_7\!:\!1,\; p_{r_1}\!:\!1,\;1^3]$};

\node[fval, anchor=south] at ([xshift=25pt,yshift=3pt]S10.north west) {$(3,4,7)$};
\node[info, anchor=north west] at ([xshift=-20pt,yshift=-4pt]S10.south)
  {$[p_2\!:\!1,\; p_6\!:\!1^2,$\\$p_7\!:\!1^2,\; p_8\!:\!1,\; p_{r_1}\!:\!2^2]$};

\node[fval, anchor=south] at ([xshift=-10pt]S12.north) {$(5,1,6)$};
\node[info, anchor=north] at ([yshift=-4pt]S12.south)
  {$[p_5\!:\!1,\; p_6\!:\!1,$\\$p_{10}\!:\!1^1,\; p_{r_1}\!:\!1,\;1^1]$};

\node[fval, anchor=south west] at ([xshift=0pt,yshift=3pt]S11.north east) {$(5,1,6)$};
\node[info, anchor=west] at ([xshift=0pt]S11.east)
  {$[p_7\!:\!1,\; p_8\!:\!1,$\\$p_9\!:\!1^1,\; p_{r_1}\!:\!1,\;1^1]$};

\node[fval, anchor=west] at ([xshift=0pt,yshift=15pt]S13.east) {$(5,1,6)$};
\node[info, anchor=north west] at ([xshift=0pt,yshift=10pt]S13.east)
  {$[p_9\!:\!1^1,\;p_{10}\!:\!1^1,$\\$p_{r_1}\!:\!2^1]$};

\node[fval, anchor=west] at ([xshift=-25pt,yshift=10pt]S14.north west) {$(6,0,6)$};
\node[info, anchor=north west] at ([xshift=-40pt,yshift=-20pt]S14.east)
  {$[p_{11}\!:\!1,\;p_{r_1}\!:\!2]$};

\draw[arc] (S0) -- node[tlbl, above] {$0$} (S1);
\draw[arc] (S1) -- node[tlbl, above] {$a$} (S2);
\draw[arc] (S2) -- node[tlbl, above] {$c$} (S6);

\draw[arc] (S1) --
    node[tlbl, left, xshift=-2pt, pos=0.30] {$b$} (S3);
\draw[arc] (S1) --
    node[tlbl, right, xshift=2pt, pos=0.30] {$c$} (S4);

\draw[arc] (S4) to
    node[tlbl, right, xshift=2pt,pos=0.35] {$b$} (S10);

\draw[arc] (S6) --
    node[tlbl, right, xshift=2pt,pos=0.40] {$b$} (S7);

\draw[arc] (S2) --
    node[tlbl, right, pos=0.45] {$b$} (S5);

\draw[arc] (S3) --
    node[tlbl, left, xshift=-2pt,pos=0.40] {$c$} (S8);
\draw[arc] (S3) --
    node[tlbl, right, xshift=2pt,pos=0.40] {$a$} (S9);

\draw[arc] (S7) --
    node[tlbl, above, xshift=4pt,pos=0.45] {$d$} (S11);
\draw[arc] (S7) --
    node[tlbl, left,xshift=-2pt, pos=0.45] {$e$} (S12);

\draw[arc] (S11) --
    node[tlbl, right, pos=0.50] {$e$} (S13);
\draw[arc] (S13) --
    node[tlbl, above, xshift=4pt, pos=0.50] {$n$} (S14);

\draw[dasharc] (S4) to[out=20, in=-145]
    node[tlbl, above, pos=0.42] {$a$} (S6);

\draw[dasharc] (S12.north east) to[out=40,in=225]
    node[tlbl, below, pos=0.5] {$d$}
    (S13.south west);

\end{tikzpicture}}
    \caption{Partial timed reachability graph of the $TTPNR$ for the running example.} \label{fig_A_star_on_RG}
    \vspace{-.3cm}\caption*{\scriptsize 
    \textit{Note:} Each state $s$ is accompanied by its timed marking $M_s$ and the tuple $(g(s),h(s),\feval(s))$, where $g(s)$ is the current system time, $h(s)$ is an admissible lower bound on the remaining time to project completion, and $\feval(s)=g(s)+h(s)$ is the resulting lower bound on the final project duration used by $A^*$. Each directed edge is labeled by the activity whose transition is fired. Gray states are generated and expanded by $A^*$, whereas white states are generated but not expanded. The double-bordered state $s_{14}$ is a goal state, in which all activities have completed. The dashed edges indicate transitions to previously generated states, illustrating duplicate-state detection. State indices reflect the generation order in the illustrative $A^*$ run. For Section~\ref{sec:RG_TTPNR}, the relevant elements are the timed markings and transition firings. The $A^*$ evaluation and expansion annotations are discussed in Section~\ref{sec:A_star}.}
\end{figure}

Let us illustrate the notation using the $TTPNR$ in Figure~\ref{fig:toyTTPNR}. At the initial node $s_0$, the timed marking is $M_{s_0}=[p_1:1,\ p_{r_1}:2]$, and the current system time is $g(s_0)=0$. This marking indicates that there is one immediately available control-flow token in place $p_1$ and two immediately available resource tokens in place $p_{r_1}$. After $t_0$ is fired, the system time remains $0$, since $\tau(t_0)=0$, and the timed marking becomes $M_{s_1}=[p_2:1,\ p_3:1,\ p_4:1,\ p_{r_1}:2]$. The input places of activities $a$, $b$, and $c$ contain the required control-flow tokens, and the resource place contains two immediately available resource tokens. Hence, each of $t_a$, $t_b$, and $t_c$ can be selected as the next scheduling decision. Suppose, for example, that $t_a$ is fired from $s_1$. This produces state $s_2$ with timed marking $M_{s_2}=[p_3:1,\ p_4:1,\ p_5:1^3,\ p_{r_1}:1,1^3]$. Transition $t_a$ consumes the token in $p_2$ and one resource token. A control-flow token becomes available in $p_5$ after $\tau(t_a)=3$ time units, and the consumed resource token returns to $p_{r_1}$ after the same delay. Although the control-flow token required by activity $b$ is immediately available in $p_3$, only one resource token is immediately available. Therefore, selecting $b$ next advances the system time to time $3$, when the resource token used by activity $a$ returns and two resource tokens are again available. Alternatively, if transition $t_b$ is fired immediately after $t_0$, the resulting state $s_3$ has timed marking $M_{s_3}=[p_2:1,\ p_4:1,\ p_6:1^2,\ p_7:1^2,\ p_{r_1}:2^2]$. Activity $b$ is fully enabled at time $0$ and has duration $\tau(t_b)=2$. Its output control-flow tokens and the two resource tokens it consumes therefore become available after two time units.

Figure~\ref{fig_A_star_on_RG} illustrates the timed markings associated with the states of the reachability graph and the way these markings evolve under transition firings. The figure also contains the evaluation annotations and node-status shading, which are explicitly explained in the $A^*$ illustration of Section~\ref{sec:A_star}.

\subsection{Mapping between Representations}\label{sec:equivalences}
The three representations, MILP, state-action, and Petri net, describe the same underlying RCPSP from different viewpoints. The MILP formulation encodes complete schedules, whereas the state-action and Petri-net representations encode intermediate execution states and their transitions. In particular, the state-action representation records explicitly which activities are completed, ongoing, or not yet started, while the Petri-net representation records the same information implicitly through timed control-flow and resource tokens.

\subsubsection{Mapping between State-Action and Petri Net Representations}
For the RCPSP-induced $TTPNR$ considered in this paper, the reachability graph $RG=(\Sset_{RG},\mathbb{E}_{RG})$ of Section~\ref{sec:RG_TTPNR} is already defined at \emph{decision states}. Each node is a timed state represented by a timed marking and the current system time, and each outgoing edge corresponds to choosing one activity to start next, with the clock advancing by the smallest amount needed to make the required control-flow and resource tokens available. Therefore, the Petri-net reachability graph and the state graph use the same transition granularity, differing only in how states are represented.

There is thus a natural correspondence between a state $s=\bigl(g(s),\,\V^c(s),\,\V^o(s),\,\tilde d(s)\bigr)$ and a timed marking at the same system time. Completed activities correspond to immediately available control-flow tokens at the appropriate successor places, ongoing activities and their remaining durations correspond to delayed control-flow and resource tokens with positive residual duration, and currently available resources correspond to immediately available tokens in the resource places. 

\begin{proposition}\label{prop:sa_pn}
The state graph $G=(\Sset,\mathbb{E})$ and the $TTPNR$ reachability graph
$RG=(\Sset_{RG},\mathbb{E}_{RG})$, where $\Sset_{RG}$ is the set of reachable timed states and $\mathbb{E}_{RG}$ is the set of directed edges induced by transition firings, are isomorphic. That is, there exists a one-to-one correspondence between the states of $G$ and the states of $RG$, $\phi:\Sset\to\Sset_{RG}$, such that:
\begin{enumerate}
    \item[\emph{(i)}] $\phi(s_0)$ is the initial RG node represented by the timed marking $M_{s_0}$ and system time $0$, where $M_{s_0}$ refines the ordinary initial marking $M_0$,
    \item[\emph{(ii)}] corresponding states have the same system time,
    \item[\emph{(iii)}] for every activity $j\in\V$, there is an edge $s \xrightarrow{\,j\,} s'$ in $G$ if and only if there is an edge $\phi(s) \xrightarrow{\,t_j\,} \phi(s')$ in $RG$ with label $L(t_j)=j$.
\end{enumerate}
\end{proposition}
\begin{proof}
The proof is presented in Appendix~\ref{app:proof_prop1}.
\end{proof}

\subsubsection{Mapping between MILP and State-Action  Representations}
A partial schedule can be mapped to a state in the state graph induced by the state-action representation, but this mapping is generally many-to-one. Different partial schedules may lead to the same state. For example, in the running example, scheduling the independent activities $a$ and $c$ in either order may produce identical states. The following proposition formalizes this relation.

\begin{proposition}\label{prop:milp_sa}
    Given an RCPSP, let $G$ be the state graph associated with the state-action representation. Then, any partial schedule $(\J,\{S_i\}_{i\in\J})$ produced by earliest-start generation corresponds to a path in the graph $G$ leading to a state $s\in\Sset$ satisfying
    $g(s)=\max_{i\in\J} S_i$ and $\V^o(s)=\{i\in\J: S_i\leq g(s)< S_i+d_i\}$, $\tilde{d}_i(s)=S_i+d_i-g(s)>0$ for all $i\in \V^o(s)$, and $\V^c(s)=\J\setminus \V^o(s)$.
    Moreover, any state $s\in \Sset$ corresponds to at least one partial schedule $(\J,\{S_i\}_{i\in\J})$ such that $\J=\V^o(s)\cup\V^c(s)$, and $S_i=g(s)+\tilde{d}_i(s)-d_i$ for all $i\in \V^o(s)$ and $S_i\leq g(s)-d_i$ for all $i\in \V^c(s)$.
\end{proposition}
\begin{proof}
The proof is presented in Appendix~\ref{app:proof_prop2}.
\end{proof}

In particular, the above proposition implies that a complete path from $s_0$ to a goal state  $s_f$ in the graph $G$ defines a schedule $\{S_i\}_{i\in\V}$ by recording the time at which each activity's transition fires. Conversely, given a feasible MILP solution with start times $\{S'_i\}_{i\in\V}$, the corresponding path in $G$ can be reconstructed by processing activities in non-decreasing order of their start times.

Note that since a state in $\Sset$ is not necessarily produced by a unique partial schedule, there may be multiple feasible schedules (paths in $G$) associated with a single state $s_f\in\Sset_f$. This many-to-one relationship between scheduling decisions and states may affect solution approaches that use these two representations. 

Thus, while $A^*$ is naturally applied to the graph $G$ associated with the state-action representation, the ordinary OR-style B\&B approach is naturally applied to a MILP or partial-schedule representation. The distinction that $A^*$ searches over states, whereas B\&B searches over a tree of partial scheduling decisions is central to the comparison below.

\section{Solution Approaches}
\label{sec:approaches}

In this section, we describe the two exact solution methods B\&B (Section~\ref{sec:bb}) and $A^*$ (Section~\ref{sec:A_star}) and provide a structural comparison of the two approaches (Section~\ref{sec:structural_comparison}). We first present each method in a generic form and then specialize it to the RCPSP. Both methods use lower bounds on the best achievable project completion time to guide their search. For the baseline comparison, we use the same critical-path lower bound and best-first (BF) node-selection rule for both methods. This design controls for the common lower bound and search order while retaining the algorithm-specific mechanisms that distinguish the two approaches. In Section~\ref{sec:empirical}, we additionally examine B\&B variants with depth-first (DF) search and linear programming (LP) relaxations that produce stronger lower bounds. 

\subsection{Branch-and-Bound}
\label{sec:bb}

In the field of OR, B\&B is an exact solution framework for combinatorial optimization problems that systematically explores the solution space while avoiding exhaustive enumeration. In the context of the RCPSP, it is naturally applied to the MILP representation and corresponding partial schedules.

Within B\&B, the explored solutions can be represented by a rooted search tree, in which the root node $N_0$ corresponds to the optimization problem to be solved, i.e., its objective function and constraints. Each other node $N_j$ of the B\&B search tree is implicitly defined by the objective function of the problem to be solved and the constraints of the problem corresponding to its parent node $N_{i}$, as well as an additional attribute (usually a constraint) that every solution to the problem at $N_j$ must satisfy, whereas this need not be the case for all solutions to the problem at $N_{i}$. A path in the search tree from the root node $N_0$ to a leaf node thus corresponds to the initial problem in $N_0$, successively extended by the constraints added along the branches of this path. Consequently, the set of feasible solutions can only shrink from a node $N_i$ to its child $N_j$, while the best achievable objective value of the corresponding subproblem cannot improve. Hence, B\&B relies on partitioning the solution space (branching) and eliminating problems (nodes) that cannot lead to an optimal solution (bounding and pruning). The overall search is typically guided by strategies for selecting the next node to be evaluated.

Since the best achievable objective value of a node is hard to calculate in general, bounds are computed at the nodes of the search tree. First, a global upper bound is the objective value of a feasible solution of the problem in $N_0$, which may be improved several times during the search as new, better feasible solutions of $N_0$ are found. Second, a local lower bound of the problem in $N_i$ is computed at each node $N_i$, i.e., a value not exceeding the best-possible objective function value that any feasible solution in $N_i$ can achieve. The computation of this bound can, for example, be performed by solving a relaxed (simplified) version of the optimization problem considered in $N_i$ that can be solved efficiently. If this bound indicates that no better solution can be found within the corresponding subtree originating in $N_i$, the node is discarded (pruned, fathomed). 

\subsubsection{Generic Framework for Minimization Problems}

Consider a combinatorial optimization problem $\min_{x \in \mathcal{X}} F(x)$ with the set of feasible solutions $\mathcal{X}$ and objective function $F: \mathcal{X} \rightarrow \mathbb{R}$. Assume that $\mathcal{X}$ is too large to be enumerated explicitly, so that it is systematically explored by means of B\&B. The root node $N_0$ represents the original problem. The remaining nodes $N_i$ are associated with subsets $\mathcal{X}(N_i) \subseteq \mathcal{X}$ of feasible solutions. Hence, the branching of a node $N_i$ corresponds to partitioning $\mathcal{X}(N_i)$ into smaller disjoint subsets $\mathcal{X}(N_i) = \bigcup_{j=i+1}^k \mathcal{X}(N_j)$, where $N_{i+1}, \dots, N_k$ are the child nodes of $N_i$. This partitioning is typically induced by fixing decision variables or imposing additional constraints.

The algorithm maintains a global upper bound $\bar F$, which corresponds to the objective value of the best feasible solution of the original problem found so far. Additionally, for each node $N_i$, a lower bound $\underbar F (N_i) \leq \min_{x \in \mathcal{X}(N_i)} F(x)$ on the objective value of all solutions in $\mathcal{X}(N_i)$ is computed. The quality of these bounds affects the efficiency of the algorithm.

A node $N_i$ is pruned if one of the following conditions holds:
\begin{itemize}
    \item $\underbar F(N_i) \geq \bar F$, since no solution in $\mathcal{X}(N_i)$ can improve the current best solution.
    \item $\mathcal{X}(N_i)=\emptyset$, i.e., the subproblem at $N_i$ is infeasible.
    \item The relaxation solved at $N_i$ yields a feasible solution $x\in\mathcal{X}$. 
\end{itemize}

The order in which nodes are explored is determined by a search strategy. Common strategies include DF search and BF search. DF search completely explores a single branch before backtracking. It aims at quickly determining feasible solutions, and thus global upper bounds, to the original problem. BF search, in contrast, selects the next node to be explored based on the lower bounds of the not yet fathomed (active) nodes. 

In summary, a generic version of B\&B can be described as follows:
\begin{enumerate}
    \item Initialize the set of active nodes with the root node (corresponding to the original problem).
	\item Set $\bar F := +\infty$, or $\bar F := F(x^0)$ if a feasible solution $x^0$ is known in advance.
	\item While there are active nodes:
	\begin{enumerate}
		\item Select a node $N_i$ according to the search strategy.
		\item Compute $\underbar F(N_i)$, e.g., by considering and solving a relaxed version of the problem in $N_i$.
		\item If $\underbar F(N_i) \geq \bar F$ or if the problem in $N_i$ has no feasible solution, prune $N_i$, i.e., delete $N_i$ from the set of active nodes, and continue.
		\item Else, if the solution $x_i$ of the relaxed version of the problem in $N_i$ is feasible for the original problem, set $\bar F \leftarrow \min\{\bar F, F(x_i)\}$, potentially update the best solution to $x_i$, prune $N_i$ and continue.
		\item Otherwise, branch on $N_i$, i.e., delete it from the set of active nodes and add its children to the set of active nodes.
	\end{enumerate}
\end{enumerate}

The correctness of this generic version of B\&B follows from the fact that only nodes that cannot contain better solutions than the currently best known solution are pruned. Therefore, an optimal solution is never discarded if no optimal solution is yet known. The efficiency of the method critically depends on the strength of the bounds, the effectiveness of the branching scheme and the search strategy. In the worst case, the algorithm may still explore the entire search space. However, in practice, strong bounds and good branching rules can lead to substantial reductions in computational effort.

\subsubsection{A Branch-and-Bound Method for the RCPSP}
As the goal of this paper is to compare baseline variants of $A^*$ and B\&B, we refrain from adapting sophisticated B\&B approaches for the RCPSP (see Section~\ref{sec:lit}). Instead, we restrict ourselves to a fairly basic version which we describe in this section.

Each node $N_i$ is associated with a partial schedule $(\J,S_{\J})$ of the considered RCPSP instance with the current time $g(N_i)=\max_{j\in \J} S_j$. The state-consistent critical-path (CP) lower bound on the project duration reachable from node $N_i$, denoted by $\feval(N_i)$ (i.e., $\underline F(N_i)$ in the generic framework), corresponds to the optimal objective function value of the relaxation of problem \eqref{rcpsp_obj}--\eqref{rcpsp_binary} that neglects all resource constraints under the given partial schedule with start times $S_{\mathcal{J}}=\{S_j\}_{j\in\mathcal{J}}$: \allowdisplaybreaks
\begin{subequations}
\begingroup 
\addtolength{\jot}{-.3em}
\begin{align}
\feval(N_i) = \min_{S_j,\,j\in \V\setminus\J}\ &S_{n}\\
\text{s.t.}\quad
& S_j + d_j \le S_l 
&& \forall (j,l) \in \E\\
& S_j \geq g(N_i)
&& \forall j \in \V\setminus\J
\end{align}
\endgroup
\end{subequations}
Since this relaxation contains only precedence constraints, $\hat f(N_i)$ is obtained by a single forward critical-path pass and can be computed without an optimization solver.

In this relaxation, the start times of already scheduled activities are fixed. Equivalently, in the timed-state representation, completed activities require no further processing, ongoing activities contribute their residual durations, and not-yet-started activities may start no earlier than $g(N_i)$, subject to precedence constraints but with resource constraints relaxed. The resulting bound satisfies $\feval(N_i)=g(N_i)+h(N_i)$, where $h(N_i)$ is a lower bound on the remaining time to project completion.

If a node $N_i$ with partial schedule $(\J,S_{\J})$ cannot be pruned, we branch over all activities that can be started at time $t \geq g(N_i)$, while respecting all precedence and resource constraints. 
Hence, a child $N_k$ of $N_i$ is associated with augmenting the partial schedule with an additional activity $j$, such that all predecessors of $j$ belong to $\J$, where activity $j$ is started at the earliest feasible time no earlier than $g(N_i)$, resulting in the partial schedule $(\J\cup\{j\},S_{\J\cup\{j\}})$ (see Section~\ref{sec:problem_formal}). 
Hence, $S_j$ is given by the solution of the following mixed integer program, where the values $y_{it}$ are fixed to one if $S_i=t$ and to zero otherwise for all $i\in \J$: 
\begin{subequations} \allowdisplaybreaks \label{eq:compute_earlier_start}
\begingroup 
\addtolength{\jot}{-.3em}
\begin{align}
\min_{S_j, y_{jt}}\ &S_{j} \label{eq:bb1}\\
\text{s.t.}\quad
& \sum_{t=g(N_i)}^{H-d_j} y_{jt} = 1 &&\\
& S_j = \sum_{t=g(N_i)}^{H-d_j} t\, y_{jt}
&& \\
& S_l + d_l \le S_j 
&& \forall (l,j) \in \E, l\in \J\\
& \sum_{l \in \J\cup\{j\}} r_{lk}
\sum_{\tau=\max\{0,t-d_l+1\}}^{\min\{t,H-d_l\}} y_{l\tau}
\le R_k && \forall t \in \T,\  k \in \R\\
& y_{jt}=0 && \forall t=0,\ldots, g(N_i)-1\\
& y_{jt}\in\{0,1\},
&& \forall t=g(N_i),\dots,H-d_j \label{eq:bbn}
\end{align}
\endgroup
\end{subequations}
This ensures that $g(N_k)\geq g(N_i)$ and that the start time $S_j$ of activity $j$ takes into account both resource and precedence constraints. We note that in practice $S_j$ can be computed analytically without the need to solve \eqref{eq:bb1}--\eqref{eq:bbn}.

Scheduling all activities sequentially yields a feasible schedule with makespan $H$. We therefore initialize the upper bound as $\bar F=H$.
We apply a BF search strategy that selects an active node $N_i$ with minimum value $\feval(N_i)$. As a tie-breaker, we choose the highest value of $g(N_i)$. If both values are tied, we apply a First In, First Out (FIFO) strategy. B\&B can also employ other node-selection strategies, but this choice aligns the B\&B node-selection rule with the BF expansion rule of $A^*$, thereby controlling for differences in search order and allowing the comparison to focus on the effects of tree-based enumeration and duplicate-state detection. 

For the running example, the procedure generates the search tree shown in Figure~\ref{fig_BB_Tree}.
The initial upper bound is $\bar F=H=10$. From the root node $N_0$, the only possible branching decision is dummy activity $0$, because it is the only activity without predecessors. This produces node $N_1$. Since activity $0$ has zero duration, $\bigl(g(N_1),h(N_1),\feval(N_1)\bigr)=(0,4,4)$. At $N_1$, activities $a$, $b$, and $c$ become available, generating nodes $N_2$, $N_3$, and $N_4$, respectively. Each has evaluation value $4$. Under the above tie-breaking rules, $N_2$ is selected first, and activity $a$ is scheduled at time $0$. Branching at $N_2$ generates $N_5$, corresponding to selecting activity $b$, and $N_6$, corresponding to selecting activity $c$. Activity $b$ requires two resource units and therefore cannot start until activity $a$ completes. Hence, $\bigl(g(N_5),h(N_5),\feval(N_5)\bigr)=(3,4,7)$. At this stage, $N_5$ is not pruned because $\feval(N_5)=7<\bar F=10$. It stays active but is not selected while nodes with smaller evaluation values exist. Activity $c$, in contrast, can start in parallel with activity $a$. Consequently, $\bigl(g(N_6),h(N_6),\feval(N_6)\bigr)=(0,4,4)$. Branching at $N_3$ generates $N_7$ and $N_8$, corresponding to selecting activities $a$ and $c$, respectively; both have evaluation triple $(2,4,6)$. Branching at $N_4$ generates $N_{9}$ by selecting activity $a$ and $N_{10}$ by selecting activity $b$. Since activity $a$ can start in parallel with activity $c$, node $N_{9}$ has evaluation triple $(0,4,4)$. By the FIFO rule, $N_6$ is expanded next. Selecting activity $b$ from $N_6$ subsequently generates $N_{11}$ with evaluation triple $(3,3,6)$. The search next processes the remaining node $N_9$ with evaluation value $4$, selecting activity $b$ to produce $N_{12}$ with evaluation triple $(3,3,6)$. Expanding $N_{11}$ generates $N_{13}$ and $N_{14}$ by scheduling activities $d$ and $e$, respectively, both with $(g,h,\feval)=(5,1,6)$. Since $g(N_{13})=g(N_{14})=5>g(N_{12})=3$, the tie-breaking rule prefers $N_{13}$ and $N_{14}$ over $N_{12}$, and FIFO selects $N_{13}$, which generates $N_{15}$ by scheduling $e$. Next, $N_{14}$ is expanded, generating $N_{16}$ by scheduling $d$. Nodes $N_{15}$ and $N_{16}$ represent the same partial schedule but remain distinct nodes in the B\&B tree. By FIFO, $N_{15}$ is selected next and generates the goal node $N_{17}$ by scheduling the dummy end activity $n$. With $(g,h,\feval)=(6,0,6)$, $N_{17}$ is selected
next, updating the upper bound to $\bar F=6$. All remaining active nodes have lower bounds of at least $6$, and the search terminates with the optimal makespan of $6$.

\begin{figure}[t]
    \centering
    \begin{tikzpicture}[
        x=1cm,
        y=1cm,
        state/.style={
            circle,
            draw=black,
            thick,
            fill=white,
            minimum size=9mm,
            inner sep=0pt,
            font=\footnotesize
        },
        expandedstate/.style={
            circle,
            draw=black,
            thick,
            fill=gray!30,
            minimum size=9mm,
            inner sep=0pt,
            font=\footnotesize
        },
        goalstate/.style={
            circle,
            draw=black,
            thick,
            fill=gray!40,
            double,
            double distance=1.1pt,
            minimum size=9mm,
            inner sep=0pt,
            font=\footnotesize
        },
        arc/.style={
            -{Stealth[length=4.5pt,width=3.5pt]},
            thick,
            shorten >=1pt,
            shorten <=1pt
        },
        tlbl/.style={
            font=\footnotesize,
            fill=white,
            inner sep=1pt
        },
        fval/.style={
            font=\scriptsize,
            fill=white,
            inner sep=1pt
        }
    ]
    \node[expandedstate] (N0) at (0,0) {$N_0$};
    \node[expandedstate] (N1) at (0,-1.5) {$N_1$};
    \node[expandedstate] (N2) at (-3.8,-3.0) {$N_2$};
    \node[expandedstate] (N3) at (0,-3.0) {$N_3$};
    \node[expandedstate] (N4) at (3.8,-3.0) {$N_4$};
    \node[state]         (N5)  at (-5.2,-4.5) {$N_5$};
    \node[expandedstate] (N6)  at (-2.7,-4.5) {$N_6$};
    \node[state]         (N8)  at (-1.0,-4.5) {$N_7$};
    \node[state]         (N9)  at (1.0,-4.5) {$N_8$};
    \node[expandedstate] (N10) at (2.7,-4.5) {$N_{9}$};
    \node[state]         (N11) at (5.2,-4.5) {$N_{10}$}; 
    \node[expandedstate] (N7)  at (-2.7,-6.0) {$N_{11}$};
    \node[state] (N12) at (2.7,-6.0) {$N_{12}$};    
    \node[expandedstate] (N13) at (-3.75,-7.5) {$N_{13}$};
    \node[expandedstate] (N14) at (-1.65,-7.5) {$N_{14}$};
    \node[expandedstate] (N15) at (-3.75,-9.0) {$N_{15}$};
    \node[state]         (N16) at (-0.9,-9.0) {$N_{16}$};
    \node[goalstate]     (N17) at (-3.75,-10.5) {$N_{17}$};
    \node[fval, anchor=west]
        at ([xshift=6pt,yshift=2pt]N0.east)
        {$(0,4,4)$};
    \node[fval, anchor=west]
        at ([xshift=6pt,yshift=2pt]N1.east)
        {$(0,4,4)$};
    \node[fval, anchor=east]
        at ([xshift=-6pt,yshift=2pt]N2.west)
        {$(0,4,4)$};
    \node[fval, anchor=west]
        at ([xshift=6pt,yshift=2pt]N3.east)
        {$(0,4,4)$};
    \node[fval, anchor=west]
        at ([xshift=6pt,yshift=2pt]N4.east)
        {$(0,4,4)$};
    \node[fval, anchor=east]
        at ([xshift=-6pt,yshift=2pt]N5.west)
        {$(3,4,7)$};
    \node[fval, anchor=east]
        at ([xshift=-6pt,yshift=2pt]N6.west)
        {$(0,4,4)$};
    \node[fval, anchor=east]
        at ([xshift=-6pt,yshift=2pt]N7.west)
        {$(3,3,6)$};
    \node[fval, anchor=south]
        at ([xshift=0pt,yshift=-10pt]N8.south)
        {$(2,4,6)$};
    \node[fval, anchor=south]
        at ([xshift=0pt,yshift=-10pt]N9.south)
        {$(2,4,6)$};
    \node[fval, anchor=east]
        at ([xshift=32pt,yshift=2pt]N10.east)
        {$(0,4,4)$};
    \node[fval, anchor=west]
        at ([xshift=6pt,yshift=2pt]N11.east)
        {$(3,4,7)$};
    \node[fval, anchor=west]
        at ([xshift=6pt,yshift=2pt]N12.east)
        {$(3,3,6)$};
    \node[fval, anchor=east]
        at ([xshift=-6pt,yshift=2pt]N13.west)
        {$(5,1,6)$};
    \node[fval, anchor=west]
        at ([xshift=6pt,yshift=2pt]N14.east)
        {$(5,1,6)$};
    \node[fval, anchor=west]
        at ([xshift=6pt,yshift=2pt]N15.east)
        {$(5,1,6)$};
    \node[fval, anchor=west]
    at ([xshift=6pt,yshift=2pt]N16.east)
    {$(5,1,6)$};
    \node[fval, anchor=west]
    at ([xshift=6pt,yshift=2pt]N17.east)
    {$(6,0,6)$};
    \draw[arc]
    (N0) --
    node[tlbl, left, pos=0.50] {$0$}
    (N1);
    \draw[arc]
        (N1) --
        node[tlbl, above left, pos=0.45] {$a$}
        (N2);
    \draw[arc]
        (N1) --
        node[tlbl, left, pos=0.50] {$b$}
        (N3);
    \draw[arc]
        (N1) --
        node[tlbl, above right, pos=0.45] {$c$}
        (N4);
    \draw[arc]
        (N2) --
        node[tlbl, above left, pos=0.48] {$b$}
        (N5);
    \draw[arc]
        (N2) --
        node[tlbl, above right, pos=0.48] {$c$}
        (N6);
    \draw[arc]
        (N6) --
        node[tlbl, left, pos=0.50] {$b$}
        (N7);
    \draw[arc]
        (N3) --
        node[tlbl, above left, pos=0.48] {$a$}
        (N8);
    \draw[arc]
        (N3) --
        node[tlbl, above right, pos=0.48] {$c$}
        (N9);
    \draw[arc]
        (N4) --
        node[tlbl, above left, pos=0.48] {$a$}
        (N10);
    \draw[arc]
        (N4) --
        node[tlbl, above right, pos=0.48] {$b$}
        (N11);
    \draw[arc]
        (N10) --
        node[tlbl, left, pos=0.50] {$b$}
        (N12);
    \draw[arc]
        (N7) --
        node[tlbl, above left, pos=0.48] {$d$}
        (N13);
    \draw[arc]
        (N7) --
        node[tlbl, above right, pos=0.48] {$e$}
        (N14);
    \draw[arc]
        (N13) --
        node[tlbl, left, pos=0.50] {$e$}
        (N15);
    \draw[arc]
    (N14) --
    node[tlbl, right=2pt, pos=0.50] {$d$}
    (N16);
    \draw[arc]
    (N15) --
    node[tlbl, left, pos=0.50] {$n$}
    (N17);
    \end{tikzpicture}%
    \caption{B\&B search tree generated by a best-first procedure.    }
    \label{fig_BB_Tree}
    \vspace{-.3cm}\caption*{\scriptsize \textit{Note:} Gray nodes are generated and expanded, and white nodes are only generated. Each node shows $(g,h,\feval)$ values. Node numbers reflect the node generation order.}
\end{figure}

Our procedure is refined by two additional elements. Observe that the above implementation may create duplicate schedules due to different insertion of activities starting at the same time. To prevent such duplication, we prune nodes which insert activities starting at the same time in non-lexicographic order. In the example of Figure~\ref{fig_BB_Tree}, node $N_{9}$ would be pruned immediately as it schedules both a and c at time 0 in non-lexicographic order. Similarly, $N_{16}$ would be pruned because activities $e$ and $d$ are scheduled at the same time in non-lexicographic order. Additionally, the schedule produced by the above implementation may not be active. Specifically, if Node $N_k$ is created from node $N_i$ by scheduling activity $j$, \eqref{eq:compute_earlier_start} only enforces earliest starting time from $g(N_i)$ onward. To verify that the schedule is active, we need to check whether or not $j$ can be scheduled prior to time $g(N_i)$. Thus, if replacing $g(N_i)$ by zero in \eqref{eq:compute_earlier_start}, produces $S_j<g(N_i)$ then the node $N_k$ would be pruned, as it results in a nonactive schedule. In Figure~\ref{fig_BB_Tree}, such pruning would have occurred if, for example, we expanded node $N_5$ to schedule activity $c$, which would have resulted in a nonactive schedule.

\subsection{\texorpdfstring{$A^*$}{A*} Search}\label{sec:A_star}

$A^*$ is a well-known BF graph-search algorithm for finding a minimum-cost path from an initial state to a goal state \citep{hart1968formal}. We first describe the generic $A^*$ framework and then specialize it to the RCPSP state graph defined in Section~\ref{sec:dp}.

\subsubsection{Generic Framework}

Consider a directed state-space graph with an initial state, a set of goal states, and nonnegative transition costs. For each generated state $s$, $A^*$ uses the evaluation function $\feval(s)=g(s)+h(s)$, where $g(s)$ is the cost of the best path currently known from the initial state to $s$, and $h(s)$ estimates the minimum remaining cost from $s$ to a goal state. The heuristic $h$ is admissible if it never overestimates this remaining cost. If, in addition, $h$ is consistent, the standard graph-search version of $A^*$ guarantees that the first goal state selected for expansion lies on an optimal path.

$A^*$ maintains two sets: \textsc{Open}, containing generated states that are candidates for expansion, and \textsc{Closed}, containing states that have already been expanded. Initially, the initial state is inserted into \textsc{Open}. At each iteration, $A^*$ selects a state from \textsc{Open} with minimum $\feval$-value. Ties may be resolved using a secondary rule, such as preferring the state with the larger $g$-value.

The selected state is expanded by generating its successor states, i.e., the states reachable through its outgoing edges. A newly generated state is inserted into \textsc{Open}; if the same state has already been generated through another path, only the better path needs to be retained. The expanded state is moved to \textsc{Closed}. The procedure continues until a goal state is selected or \textsc{Open} becomes empty.

\subsubsection{\texorpdfstring{$A^*$}{A*} for the RCPSP}

We now apply this generic framework to the RCPSP state graph $G$ defined in Section~\ref{sec:dp}. Its initial state is $s_0$, its goal states are the states in $\Sset_f$, and each edge corresponds to selecting one activity to be scheduled next. The transition cost from $s$ to its successor $s'$ is $\Delta(s,s')=g(s')-g(s)$. Consequently, the accumulated path cost from $s_0$ to state $s$ equals the current system time $g(s)$.

For the RCPSP, the evaluation function is therefore $\feval(s)=g(s)+h(s)$, where $h(s)$ is an admissible lower bound on the remaining time to project completion. When a scheduling decision is made, time advances to the earliest moment at which the selected activity can begin, namely when its precedence and resource requirements are satisfied. We use the same tie-breaking rule as for B\&B: among states with equal $\feval$-values, we select the state with the largest $g$-value and remaining ties are resolved by FIFO. Furthermore, we use the same CP lower bound defined in Section~\ref{sec:bb}. Thus, $h(s)=\feval(s)-g(s)$ provides an admissible lower bound on the remaining time to project completion. We also consider the two pruning refinements introduced at the end of Section~\ref{sec:bb}: lexicographic pruning of alternative orderings of activities starting at the same time, and active-schedule pruning of scheduling decisions that could have been made earlier without delaying previously scheduled activities.

In Figure~\ref{fig_A_star_on_RG}, $A^*$ searches the reachability graph induced by the $TTPNR$ in Figure~\ref{fig:toyTTPNR}, which is equivalent to $G$ by Proposition~\ref{prop:sa_pn}. Each node represents a timed state consisting of the current system time $g(s)$ and a timed marking $M_s$ specifying which control-flow and resource tokens are immediately available and which are delayed by ongoing activities. Each outgoing arc represents one feasible next scheduling decision: choosing an activity to start next. This decision may leave $g(s)$ unchanged if the activity can start immediately, or advance time to the earliest moment at which its precedence and resource requirements are satisfied.

For every generated state $s$, $A^*$ computes $\feval(s)=g(s)+h(s)$ and selects an open state with minimum $\feval$ for expansion. From the initial state $s_0$, only dummy activity $0$ is enabled, so expanding $s_0$ generates only $s_1$. From $s_1$, activities $a$, $b$, and $c$ are all feasible next decisions, and therefore three successor states $s_2$, $s_3$, and $s_4$ are generated. In this way, the general BF principle of $A^*$ is instantiated as a search over timed project states.

To illustrate the calculation of $\feval(s)$, consider state $s_{10}$ in Figure~\ref{fig_A_star_on_RG}, whose timed marking is represented by
$M_{s_{10}}=[p_5:1,\ p_6:1^2,\ p_7:1^2,\ p_8:1,\ p_{r_1}:2^2]$. At state $s_6$, the decisions to execute activities $0$, $a$, and $c$ have already been taken, and $\bigl(g(s_6),h(s_6)\bigr)=(0,4)$. Activity $b$ is the only selectable activity with respect to the precedence constraints. It requires two resource tokens, which become available when activities $a$ and $c$ complete after three time units.  Consequently, $g(s_{10})=S_b=3$. The remaining critical-path duration through activities $\{b,d,e,n\}$ is three time units, so $h(s_{10})=3$ and $\feval(s_{10})=g(s_{10})+h(s_{10})=6$.

\subsection{Structural Comparison of Branch-and-Bound and \texorpdfstring{$A^*$}{A*}}
\label{sec:structural_comparison}

Although B\&B and $A^*$ both systematically explore alternative scheduling decisions using lower bounds to guide the search, they differ in several fundamental respects. In this section, we identify the key structural differences and discuss their implications for solving the RCPSP. The results are summarized in Table~\ref{tab:comparison}. 

The B\&B algorithm searches a tree whose nodes are generated by sequences of scheduling decisions, whereas $A^*$ searches a graph of timed project states. Every B\&B node induces a timed state containing the information required to determine the remaining feasible scheduling decisions. Different branches of the B\&B tree may induce the same timed state and therefore represent the same remaining scheduling problem. Tree-based B\&B retains such branches as separate nodes. For example, the decision sequences $\langle 0,a,c\rangle$ and $\langle 0,c,a\rangle$ generate different paths in the B\&B tree, although they may result in the same scheduled activities, current time, ongoing activities, and resource availability. Graph-based $A^*$, in contrast, identifies nodes by their timed-state representation. When a generated successor is already present in \textsc{Open} or \textsc{Closed}, it is recognized as a duplicate and is not retained as a separate search node. The $A^*$ state graph can therefore be viewed as the graph obtained from the B\&B tree by merging nodes that induce the same timed project state. Consequently, under otherwise identical search conditions, tree-based B\&B may generate and evaluate the same remaining scheduling problem several times, whereas graph-based $A^*$ represents it once. The potential redundancy increases with project concurrency, because independent activities can be selected in many different orders while leading to the same timed state.

\begin{table}[t]
\centering
\caption{Structural comparison of B\&B and $A^*$ for the RCPSP.}
\label{tab:comparison} \footnotesize \renewcommand{\arraystretch}{0.9}
\begin{tabular}{lp{5cm}p{6cm}} 
\toprule
{Aspect} & {B\&B} & {$A^*$ Search} \\
\hline
Search structure & Tree & Graph \\
Duplicate detection & No & Yes (via \textsc{Open} and \textsc{Closed} lists) \\
Upper bound & Maintained and updated & Not used \\
Lower bound & $g + h$ (same) & $g + h$ (same) \\
Pruning mechanism & Explicit: bound ($\feval(N_i) \geq \bar F$),
infeasibility, or a feasible solution & Implicit: best-first + admissibility \\
$h$ evaluations per distinct state & Possibly $>1$ & $1$ (with duplicate check before evaluation) \\
Node selection & Flexible (best-first, depth-first, \ldots) & Best-first (required) \\
\bottomrule
\end{tabular}%
\end{table}

We now turn our attention to the bounding mechanisms of the algorithms. In their RCPSP specializations, both methods use the same evaluation $\feval=g+h$, but they differ in how this bound affects the search. B\&B maintains a global upper bound $\bar F$, updated whenever a complete schedule is found, and explicitly prunes any node $N_i$ with $\feval(N_i)\geq\bar F$. The interplay between the improving upper bound and the lower bounds of individual nodes is central to its effectiveness. Our $A^*$ implementation does not maintain an upper bound. Instead, BF expansion with an admissible and consistent heuristic guarantees that the first goal state selected from \textsc{Open} is optimal. Hence, $A^*$ does not exploit an incumbent solution for pruning in the way that B\&B does.

The methods furthermore differ with respect to the heuristic computation efficiency. With state-based duplicate detection, $A^*$ retains a single search record for each distinct state, and with a consistent heuristic, each state needs to be expanded at most once. If duplicate detection is performed before evaluating $h$, the heuristic is therefore computed only once per distinct retained state. In B\&B, distinct tree nodes can induce the same timed state, so $h$ may be evaluated multiple times for the same remaining scheduling problem.

B\&B permits flexible node-selection strategies such as depth-first, breadth-first, or hybrid. This choice interacts with the upper-bound update mechanism. $A^*$, on the other hand, requires BF selection for its optimality guarantee. When both use BF selection with the same heuristic, their expansion orders are closely related but not identical since $A^*$ may expand fewer nodes due to duplicate detection.

With respect to the memory requirements, $A^*$ stores all expanded states in a \textsc{Closed} list and all frontier states in an \textsc{Open} list. DF B\&B avoids most of this cost because it only keeps the children of the nodes on the current path from the root in memory. 
However, when B\&B uses BF selection, it has to maintain all unexpanded nodes, resulting in memory requirements comparable to $A^*$. In this setting, the \textsc{Closed} list is an additional cost of $A^*$ that can also save memory by preventing duplicate states from entering the \textsc{Open} list.

\section{Empirical Study} \label{sec:empirical}
The purpose of this section is to examine the differences between the ordinary B\&B and \(A^*\) algorithms presented above across a broad range of instances. To this end, we consider the 480 instances in the J30 benchmark from PSPLIB \citep{kolisch1997psplib}, which consists of RCPSP instances with 30 activities and varying precedence structures and resource requirements. Note that we do not aim to establish that the algorithms are competitive with state-of-the-art RCPSP solvers. 

All methods were implemented in Python 3.13.13, and experiments were conducted on a single core of an Intel Xeon Gold 6354 @ 3.00 GHz, without parallelization, with a time limit of 2 hours and a memory limit of 64 GB. All code, data, and results are available at
\url{https://github.com/ShimShtern/rcpsp-astar-vs-branch-and-bound}
(version~v1.0).

\begin{figure}[t]
    \centering
    \includegraphics[width=\linewidth]{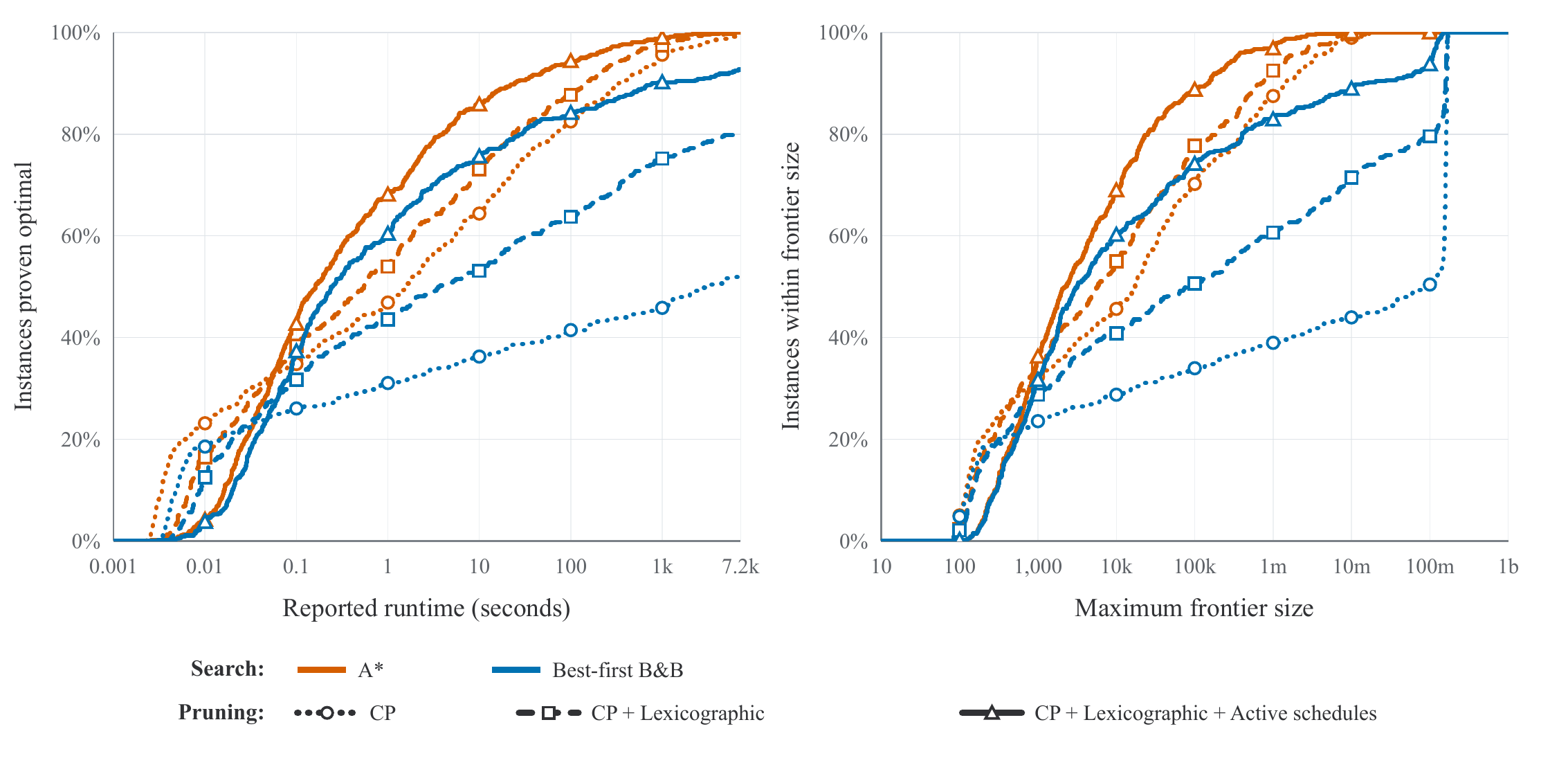}
    \caption{Time and memory behavior of $A^*$ and best-first B\&B under different pruning configurations.}
    \label{fig:time_memory_comparison_basic}
\end{figure}

We start by comparing $A^*$ and B\&B using the same CP lower bound described in  Section~\ref{sec:bb}. Thus, the only difference between BF B\&B and \(A^*\) is duplicate identification: \(A^*\) abstracts partial schedules into states and tracks duplicate states. Figure~\ref{fig:time_memory_comparison_basic} shows the effect of applying the lexicographic ordering and active-schedule pruning refinements discussed in Section~\ref{sec:bb} to both methods. The left panel depicts the number of instances proven optimal as a function of runtime, and the right panel depicts the number of instances whose maximum frontier size, i.e., the largest number of open nodes held simultaneously during the search, does not exceed a given value, which serves as a proxy for memory consumption. 

\begin{figure}[hpbt]
    \centering
    \includegraphics[width=\textwidth]{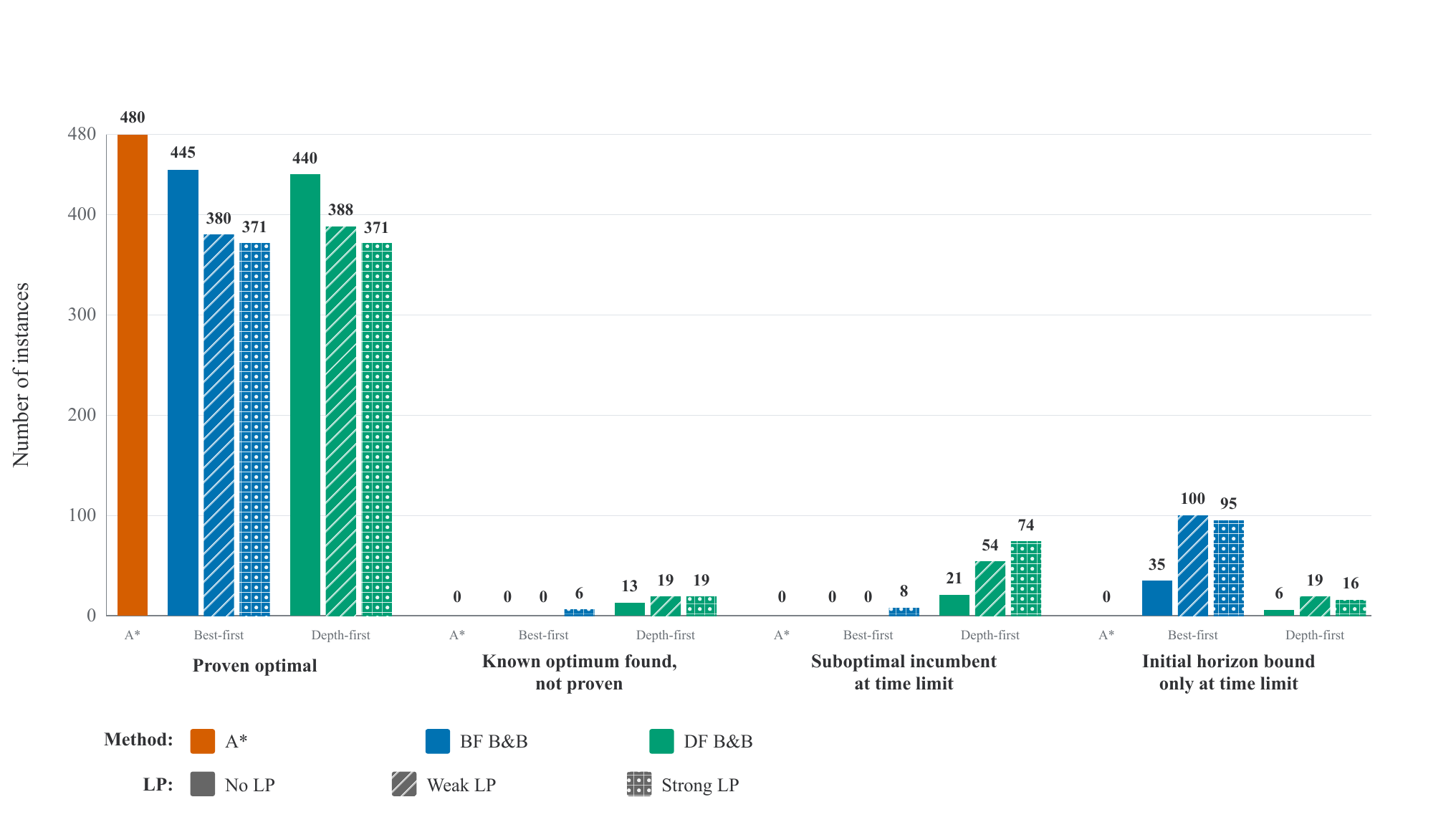}
    \caption{Status counts for each method.}
    \label{fig:bb_status_count}
\end{figure}

Without incorporating the pruning rules, BF B\&B results in 230 instances terminated upon reaching the memory limit and one upon reaching the time limit. With $A^*$, in contrast, only four instances did not reach optimality, all due to the time limit. This indicates substantially higher memory consumption for B\&B, partly due to equivalent partial schedules appearing in different parts of the search tree. Adding lexicographic pruning aided both methods. For $A^*$, three of the four previously time-limited instances reached an optimal solution. For B\&B, the number of instances terminating due to the memory limit decreased to 82, and the number terminating due to the time limit increased to 11. Lexicographic pruning therefore removes repeated occurrences of the same partial schedule across branches, but not distinct partial schedules corresponding to the same state. Among the instances reaching optimality after applying lexicographic pruning, the median runtime reduction was $93.3\%$ for B\&B, compared with $32.1\%$ for $A^*$. Active-schedule pruning further improved both methods. For $A^*$, the median reduction in runtime to optimality was $52.8\%$ and all instances were solved to optimality. For B\&B, it eliminated all memory-limit terminations, although the number of instances reaching the time limit increased to 35, and yielded a median runtime reduction of $82.4\%$. Because $A^*$ already reduces the search space through duplicate-state identification, these pruning mechanisms yield smaller additional gains than for B\&B. The right panel of Figure~\ref{fig:time_memory_comparison_basic} shows that, regardless of the pruning rules applied, $A^*$ requires a frontier that is no larger than that of BF B\&B on almost all instances. The differences become less pronounced as the pruning rules are introduced, since they reduce the frontier of B\&B considerably more than that of $A^*$. Specifically, the median ratio between the maximum frontier sizes of BF B\&B and $A^*$ decreases from about 140 without pruning rules to about 10 with lexicographic pruning and to about 1.5 with both pruning rules. Nevertheless, on the hardest instances the gap remains large, as the largest frontier of BF B\&B with both pruning rules is still more than 20 times that of $A^*$.

A potential advantage of B\&B over \(A^*\) is its ability to prune using incumbent upper bounds. We therefore examine whether exploiting this mechanism can improve the relative performance of B\&B. More specifically, we compare $A^*$ to B\&B with DF search, where the child with the better lower bound is explored first. DF search may obtain an incumbent before reaching the optimum. Additionally, we utilize the partial schedule to obtain weak and strong LP relaxations, both solved using Gurobi 13.0.3. The weak LP does not add constraints to the integer formulation, whereas the strong LP models each precedence relation through one constraint per time period, resulting in tighter lower bounds. The two formulations are detailed in Appendix~\ref{app:lp}. They provide stronger lower bounds for B\&B and, when feasible, incumbent upper bounds. The status counts for the different methods are summarized in Figure~\ref{fig:bb_status_count}. 
Without LP relaxations, DF B\&B finds suboptimal incumbent solutions for 21 instances and optimal solutions without proof of optimality for 13 additional instances out of the 40 instances terminating without proving optimality. Thus, although DF search is less effective at proving optimality, it provides useful incumbent solutions for many instances that remain unsolved within the time limit.
By construction, BF B\&B without LP relaxations produces no incumbents before the optimum, except for the initial one. In our experiments, the weak LP also yielded none, while the strong LP yielded fewer than DF B\&B. Additional results, including a breakdown of unresolved instances by PSPLIB parameter subclass and a comparison of the final upper and lower bounds of the B\&B variants, are provided in Appendix~\ref{app:additional}.

\begin{figure}[pth]
    \centering
    \includegraphics[width=\linewidth]{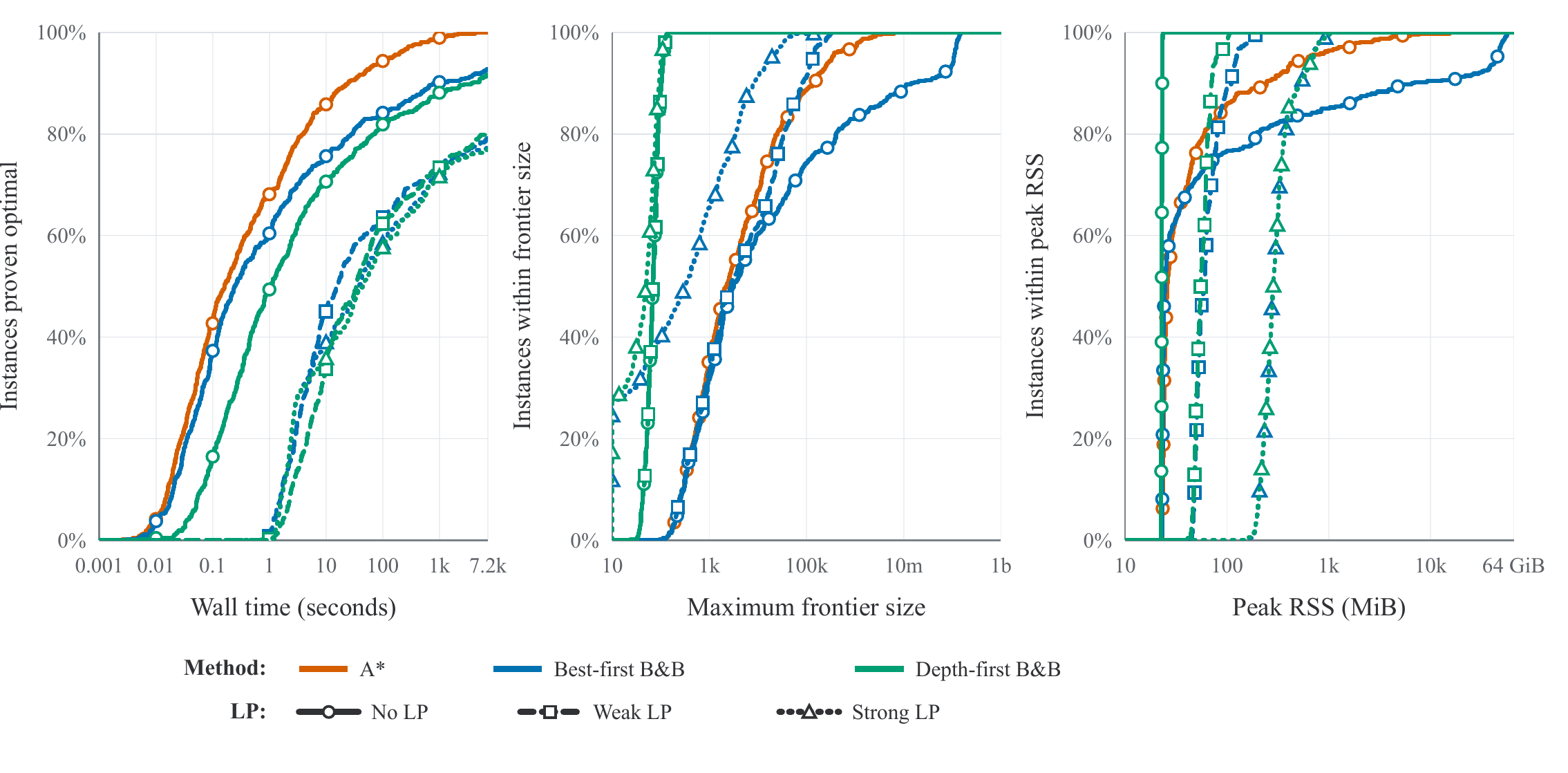}
    \caption{Time and memory performance comparing $A^*$ with BF and DF B\&B, with and without LP relaxations.}    \label{fig:time_memory_comparison_advanced}
\end{figure}

Figure~\ref{fig:time_memory_comparison_advanced} shows the time and memory profiles of the different methods over all instances, and Table~\ref{tab:astar-bb-computational-summary} summarizes the performance on the 365 instances for which all methods proved optimality.
The profiles in Figure~\ref{fig:time_memory_comparison_advanced} show that the DF search is generally much more memory efficient but slower than the BF search. Somewhat surprisingly, they also show that both LP variants reduce overall performance: although they strengthen the bounds, their computational overhead substantially increases runtime and causes more instances to terminate without proving optimality.
Table~\ref{tab:astar-bb-computational-summary} provides a closer look at the instances solved by all methods. It confirms the memory advantage of the DF search, whose maximum frontier is much smaller than that of the BF search. It further shows that the strong LP explores the tree much more efficiently than the weak LP and the no-LP variant, as reflected in the number of generated nodes, but incurs substantially higher runtime and memory usage.

\begin{table*}[t]
  \centering
  \caption{Computational summary on the 365 PSPLIB J30 instances for which
  every method proved optimality. }
  \label{tab:astar-bb-computational-summary}
  \footnotesize
  \begin{tabular}{@{}lllrrrrr@{}}
    \toprule
    Method &
    Search &
    \shortstack{LP\\relaxation} &
    \shortstack{Wall time (s)\\median {[IQR]}} &
    \shortstack{Expanded\\median {[IQR]}} &
    \shortstack{Generated\\median {[IQR]}} &
    \shortstack{Max. frontier\\median {[IQR]}} &
    \shortstack{Peak RSS (MiB)\\median {[IQR]}} \\
    \midrule
    A$^*$
      & --
      & --
      & \shortstack[r]{0.074\\{[0.026--0.300]}}
      & \shortstack[r]{1,147\\{[465--3,967]}}
      & \shortstack[r]{2,353\\{[922--8,076]}}
      & \shortstack[r]{1,111\\{[463--4,070]}}
      & \shortstack[r]{24.7\\{[23.8--28.2]}} \\
    B\&B
      & BF
      & No
      & \shortstack[r]{0.111\\{[0.033--0.560]}}
      & \shortstack[r]{1,330\\{[472--6,622]}}
      & \shortstack[r]{2,982\\{[986--13,263]}}
      & \shortstack[r]{1,548\\{[498--6,493]}}
      & \shortstack[r]{23.8\\{[23.4--26.0]}} \\
      & BF
      & Weak
      & \shortstack[r]{6.831\\{[2.742--28.366]}}
      & \shortstack[r]{1,147\\{[430--4,742]}}
      & \shortstack[r]{2,526\\{[887--11,158]}}
      & \shortstack[r]{1,315\\{[477--5,418]}}
      & \shortstack[r]{53.7\\{[49.7--62.0]}} \\
      & BF
      & Strong
      & \shortstack[r]{8.963\\{[2.413--79.996]}}
      & \shortstack[r]{69\\{[0--951]}}
      & \shortstack[r]{207\\{[1--2,445]}}
      & \shortstack[r]{81\\{[1--672]}}
      & \shortstack[r]{270.8\\{[231.4--312.8]}} \\
      & DF
      & No
      & \shortstack[r]{0.433\\{[0.118--2.215]}}
      & \shortstack[r]{6,042\\{[1,985--29,354]}}
      & \shortstack[r]{7,476\\{[2,262--40,583]}}
      & \shortstack[r]{63\\{[51--76]}}
      & \shortstack[r]{22.9\\{[22.7--23.0]}} \\
      & DF
      & Weak
      & \shortstack[r]{11.570\\{[4.554--59.114]}}
      & \shortstack[r]{5,295\\{[1,748--19,337]}}
      & \shortstack[r]{6,741\\{[1,905--31,502]}}
      & \shortstack[r]{64\\{[51--79]}}
      & \shortstack[r]{52.4\\{[49.0--60.1]}} \\
      & DF
      & Strong
      & \shortstack[r]{11.727\\{[2.385--86.468]}}
      & \shortstack[r]{224\\{[0--1,905]}}
      & \shortstack[r]{420\\{[1--3,808]}}
      & \shortstack[r]{32\\{[1--58]}}
      & \shortstack[r]{271.0\\{[230.5--314.6]}} \\
    \bottomrule
  \end{tabular}
  \vspace{-.1cm}\caption*{\scriptsize \textit{Note:} Every column uses the same common subset of 365 instances. Values in brackets are interquartile ranges. Wall time includes LP overhead and measures time until optimality was proven. Max. frontier is the largest candidate-node list observed during a run. DF denotes depth-first search with local best-lower-bound child ordering.}
\end{table*}

Overall, $A^*$ outperforms the B\&B variants. Its state-space representation and duplicate elimination reduce the effective search space enough to outweigh the additional upper- and lower-bound mechanisms available to B\&B. 

\section {Conclusions and Future Research} \label{sec:summary}

This paper investigates the computational differences between B\&B and $A^*$ using the RCPSP as a common problem. We establish the relationships between MILP, state-action, and timed Petri-net representations and compare basic implementations of both methods. By using the same lower bound and BF node-selection rule, we examine the effect of their different search mechanisms. In particular, tree-based B\&B may repeatedly explore the same remaining scheduling problem through different decision sequences, whereas graph-based $A^*$ identifies equivalent states and avoids their repeated exploration.

Experiments on 480 PSPLIB J30 instances demonstrate the importance of this difference. With lexicographic and active-schedule pruning, $A^*$ solved all 480 instances to optimality, compared with 445 for BF B\&B under the same time and memory limits. Although these pruning rules substantially improved B\&B, they did not eliminate the performance gap. DF B\&B required less memory and found feasible schedules for many instances where the search could not be completed within the time limit. Stronger LP relaxations reduced the number of generated nodes but increased computational time, illustrating the trade-off between bound
strength and evaluation effort.

These results highlight the importance of how the search space is represented and explored. They also illustrate different advantages of the two approaches: duplicate-state detection can considerably reduce redundant exploration, while DF search and incumbent-based pruning offer opportunities to control memory use and obtain feasible solutions early. Such insights can help researchers select appropriate search methods according to the characteristics and requirements of the problem. However, our findings concern the implementations and benchmark instances considered here and do not claim the general superiority of either method.

Several directions for future research arise from this study. Since both approaches admit numerous algorithmic enhancements, an important direction is to investigate the transfer of effective mechanisms between them. Examples include introducing duplicate-state detection into B\&B, incorporating incumbent-based pruning into $A^*$, applying constraint propagation to restrict scheduling decisions and strengthen bounds, and using lazy heuristic evaluation to reduce unnecessary computation.
Finally, extending the analysis to other combinatorial optimization problems would help quantify the findings across additional problem domains.

\section*{Acknowledgments}
Funded by the Deutsche Forschungsgemeinschaft -- Project number 524165670.

\section*{Use of AI Tools}
During the preparation of this work the authors used ChatGPT and Claude to assist with coding, writing, and editing. After using these tools, the authors reviewed and edited the content as needed and take full responsibility for the content of the article.

\bibliographystyle{plainnat}
\bibliography{references}

\begin{thebibliography}{37}
\providecommand{\natexlab}[1]{#1}
\providecommand{\url}[1]{\texttt{#1}}
\expandafter\ifx\csname urlstyle\endcsname\relax
  \providecommand{\doi}[1]{doi: #1}\else
  \providecommand{\doi}{doi: \begingroup \urlstyle{rm}\Url}\fi

\bibitem[Artigues et~al.(2026)Artigues, Hartmann, and Vanhoucke]{artigues2026fifty}
Christian Artigues, S{\"o}nke Hartmann, and Mario Vanhoucke.
\newblock Fifty years of research on resource-constrained project scheduling explored from different perspectives.
\newblock \emph{European Journal of Operational Research}, 328\penalty0 (2):\penalty0 367--389, 2026.
\newblock \doi{10.1016/j.ejor.2025.03.024}.

\bibitem[Bell and Park(1990)]{bell1990solving}
Colin~E. Bell and Kwangho Park.
\newblock Solving resource-constrained project scheduling problems by {$A^*$} search.
\newblock \emph{Naval Research Logistics}, 37\penalty0 (1):\penalty0 61--84, 1990.
\newblock \doi{10.1002/1520-6750(199002)37:1<61::AID-NAV3220370104>3.0.CO;2-S}.

\bibitem[Blazewicz et~al.(1983)Blazewicz, Lenstra, and Kan~Rinnooy]{blazewicz1983scheduling}
Jacek Blazewicz, Jan~Karel Lenstra, and A.~H.~G. Kan~Rinnooy.
\newblock Scheduling subject to resource constraints: classification and complexity.
\newblock \emph{Discrete Applied Mathematics}, 5\penalty0 (1):\penalty0 11--24, 1983.
\newblock \doi{10.1016/0166-218X(83)90012-4}.

\bibitem[Brucker and Knust(2003)]{brucker2003lower}
Peter Brucker and Sigrid Knust.
\newblock Lower bounds for resource-constrained project scheduling problems.
\newblock \emph{European Journal of Operational Research}, 149\penalty0 (2):\penalty0 302--313, 2003.
\newblock \doi{10.1016/S0377-2217(02)00762-2}.

\bibitem[Brucker et~al.(1998)Brucker, Knust, Schoo, and Thiele]{brucker1998branch}
Peter Brucker, Sigrid Knust, Arno Schoo, and Olaf Thiele.
\newblock A branch and bound algorithm for the resource-constrained project scheduling problem.
\newblock \emph{European Journal of Operational Research}, 107\penalty0 (2):\penalty0 272--288, 1998.
\newblock \doi{10.1016/S0377-2217(97)00335-4}.

\bibitem[Brucker et~al.(1999)Brucker, Drexl, M{\"o}hring, Neumann, and Pesch]{brucker1999resource}
Peter Brucker, Andreas Drexl, Rolf M{\"o}hring, Klaus Neumann, and Erwin Pesch.
\newblock Resource-constrained project scheduling: Notation, classification, models, and methods.
\newblock \emph{European Journal of Operational Research}, 112\penalty0 (1):\penalty0 3--41, 1999.
\newblock \doi{10.1016/S0377-2217(98)00204-5}.

\bibitem[Chen et~al.(2024)Chen, W., Tian, and Ma]{ch24rcp}
G.~Chen, He; W., Y.~Tian, and K.~Ma.
\newblock Resource-constrained project scheduling with multiple states: Bi-objective optimization model and case study of aircraft maintenance.
\newblock \emph{Computers \& Industrial Engineering}, 191:\penalty0 110169, 2024.
\newblock \doi{10.1016/j.cie.2024.110169}.

\bibitem[Coelho and Vanhoucke(2018)]{coelho2018exact}
Jos{\'e} Coelho and Mario Vanhoucke.
\newblock An exact composite lower bound strategy for the resource-constrained project scheduling problem.
\newblock \emph{Computers \& Operations Research}, 93:\penalty0 135--150, 2018.
\newblock \doi{10.1016/j.cor.2018.01.017}.

\bibitem[Cohen(2024)]{cohen2024planning}
Izack Cohen.
\newblock Data-driven project planning: An integrated network learning, process mining, and constraint relaxation approach in favor of scheduling recurring projects.
\newblock \emph{IEEE Transactions on Engineering Management}, 71:\penalty0 7719--7729, 2024.
\newblock \doi{10.1109/TEM.2024.3382727}.

\bibitem[Demeulemeester and Herroelen(1992)]{demeulemeester1992branch}
Erik Demeulemeester and Willy Herroelen.
\newblock A branch-and-bound procedure for the multiple resource-constrained project scheduling problem.
\newblock \emph{Management Science}, 38\penalty0 (12):\penalty0 1803--1818, 1992.
\newblock \doi{10.1287/mnsc.38.12.1803}.

\bibitem[Dorndorf et~al.(2000)Dorndorf, Pesch, and Phan-Huy]{dorndorf2000branch}
Ulrich Dorndorf, Erwin Pesch, and To{\`a}n Phan-Huy.
\newblock A branch-and-bound algorithm for the resource-constrained project scheduling problem.
\newblock \emph{Mathematical Methods of Operations Research}, 52\penalty0 (3):\penalty0 413--439, 2000.
\newblock \doi{10.1007/s001860000091}.

\bibitem[Guo et~al.(2023)Guo, Vanhoucke, and Coelho]{guo2023prediction}
Weikang Guo, Mario Vanhoucke, and Jos{\'e} Coelho.
\newblock A prediction model for ranking branch-and-bound procedures for the resource-constrained project scheduling problem.
\newblock \emph{European Journal of Operational Research}, 306\penalty0 (2):\penalty0 579--595, 2023.
\newblock \doi{10.1016/j.ejor.2022.08.042}.

\bibitem[Hart et~al.(1968)Hart, Nilsson, and Raphael]{hart1968formal}
Peter~E. Hart, Nils~J. Nilsson, and Bertram Raphael.
\newblock A formal basis for the heuristic determination of minimum cost paths.
\newblock \emph{IEEE Transactions on Systems Science and Cybernetics}, 4\penalty0 (2):\penalty0 100--107, 1968.
\newblock \doi{10.1109/TSSC.1968.300136}.

\bibitem[Hillion and Proth(1989)]{hillion1989petri}
H.~P. Hillion and J.-M. Proth.
\newblock Using timed petri nets for the scheduling of job-shop systems.
\newblock \emph{Engineering Costs and Production Economics}, 17\penalty0 (1-4):\penalty0 149--154, 1989.
\newblock \doi{10.1016/0167-188X(89)90064-5}.

\bibitem[Huang et~al.(2023)Huang, Zhou, Lu, and Abusorrah]{huang2023survey}
Bo~Huang, Mengchu Zhou, Xiaoyu~Sean Lu, and Abdullah Abusorrah.
\newblock Scheduling of resource allocation systems with timed petri nets: A survey.
\newblock \emph{ACM Computing Surveys}, 55\penalty0 (11):\penalty0 1--27, 2023.
\newblock \doi{10.1145/3570326}.

\bibitem[Ibaraki(1978)]{ibaraki1978branch}
Toshihide Ibaraki.
\newblock Branch-and-bound procedure and state-space representation of combinatorial optimization problems.
\newblock \emph{Information and Control}, 36\penalty0 (1):\penalty0 1--27, 1978.
\newblock \doi{10.1016/S0019-9958(78)90197-3}.

\bibitem[Kolisch(1996)]{kolisch1996serial}
Rainer Kolisch.
\newblock Serial and parallel resource-constrained project scheduling methods revisited: Theory and computation.
\newblock \emph{European Journal of Operational Research}, 90\penalty0 (2):\penalty0 320--333, 1996.
\newblock \doi{10.1016/0377-2217(95)00357-6}.

\bibitem[Kolisch and Sprecher(1997)]{kolisch1997psplib}
Rainer Kolisch and Arno Sprecher.
\newblock {PSPLIB} - {A} project scheduling problem library: {OR Software - ORSEP Operations Research Software Exchange Program}.
\newblock \emph{European Journal of Operational Research}, 96\penalty0 (1):\penalty0 205--216, 1997.
\newblock \doi{10.1016/S0377-2217(96)00170-1}.

\bibitem[Kon{\'e} et~al.(2011)Kon{\'e}, Artigues, Lopez, and Mongeau]{kone2011event}
Oumar Kon{\'e}, Christian Artigues, Pierre Lopez, and Marcel Mongeau.
\newblock Event-based {MILP} models for resource-constrained project scheduling problems.
\newblock \emph{Computers \& Operations Research}, 38\penalty0 (1):\penalty0 3--13, 2011.
\newblock \doi{10.1016/j.cor.2009.12.011}.

\bibitem[Kumar and Kanal(1983)]{kumar1983general}
Vipin Kumar and Laveen~N. Kanal.
\newblock A general branch and bound formulation for understanding and synthesizing and/or tree search procedures.
\newblock \emph{Artificial Intelligence}, 21\penalty0 (1--2):\penalty0 179--198, 1983.
\newblock \doi{10.1016/S0004-3702(83)80009-5}.

\bibitem[Labat and Pomerol(2003)]{labat2003branch}
Jean-Marc Labat and Jean-Charles Pomerol.
\newblock Are branch and bound and {$A^*$} algorithms identical?
\newblock \emph{Journal of Heuristics}, 9\penalty0 (2):\penalty0 131--143, 2003.
\newblock \doi{10.1023/A:1022573412940}.

\bibitem[Lee and DiCesare(1994)]{lee1994scheduling}
Doo~Yong Lee and Frank DiCesare.
\newblock Scheduling flexible manufacturing systems using petri nets and heuristic search.
\newblock \emph{IEEE Transactions on Robotics and Automation}, 10\penalty0 (2):\penalty0 123--132, 1994.
\newblock \doi{10.1109/70.282537}.

\bibitem[Li et~al.(2024)Li, Zhu, Zheng, and Xie]{li24sps}
H.~Li, H.~Zhu, L.~Zheng, and F.~Xie.
\newblock Software project scheduling under activity duration uncertainty.
\newblock \emph{Annals of Operations Research}, 338\penalty0 (1):\penalty0 477--512, 2024.
\newblock \doi{10.1007/s10479-023-05343-0}.

\bibitem[Lublin et~al.(2026)Lublin, Atzmon, and Cohen]{lublin2026petri}
Ido Lublin, Dor Atzmon, and Izack Cohen.
\newblock Petri net induced heuristic search for resource constrained scheduling.
\newblock \emph{Proceedings of the International Symposium on Combinatorial Search}, 19\penalty0 (1):\penalty0 234--239, 2026.
\newblock \doi{10.1609/socs.v19i1.43094}.

\bibitem[Mej{\'i}a and Odrey(2005)]{mejia2005petri}
Gonzalo Mej{\'i}a and Nicholas~G. Odrey.
\newblock An approach using {P}etri nets and improved heuristic search for manufacturing system scheduling.
\newblock \emph{Journal of Manufacturing Systems}, 24\penalty0 (2):\penalty0 79--92, 2005.
\newblock \doi{10.1016/S0278-6125(05)80009-3}.

\bibitem[Morin and Marsten(1976)]{morin1976branch}
Thomas~L. Morin and Roy~E. Marsten.
\newblock Branch-and-bound strategies for dynamic programming.
\newblock \emph{Operations Research}, 24\penalty0 (4):\penalty0 611--627, 1976.
\newblock \doi{10.1287/opre.24.4.611}.

\bibitem[Murata(1989)]{murata1989petri}
Tadao Murata.
\newblock Petri nets: Properties, analysis and applications.
\newblock \emph{Proceedings of the IEEE}, 77\penalty0 (4):\penalty0 541--580, 1989.
\newblock \doi{10.1109/5.24143}.

\bibitem[Nau et~al.(1984)Nau, Kumar, and Kanal]{nau1984general}
Dana~S. Nau, Vipin Kumar, and Laveen~N. Kanal.
\newblock General branch and bound, and its relation to {$A^*$} and {$AO^*$}.
\newblock \emph{Artificial Intelligence}, 23\penalty0 (1):\penalty0 29--58, 1984.
\newblock \doi{10.1016/0004-3702(84)90004-3}.

\bibitem[Nazareth et~al.(1999)Nazareth, Verma, Bhattacharya, and Bagchi]{nazareth1999breadth}
Terence Nazareth, Sanjay Verma, Subir Bhattacharya, and Amitava Bagchi.
\newblock The multiple resource constrained project scheduling problem: A breadth-first approach.
\newblock \emph{European Journal of Operational Research}, 112\penalty0 (2):\penalty0 347--366, 1999.
\newblock \doi{10.1016/S0377-2217(97)00402-5}.

\bibitem[Patterson(1984)]{patterson1984comparison}
James~H. Patterson.
\newblock A comparison of exact approaches for solving the multiple constrained resource, project scheduling problem.
\newblock \emph{Management Science}, 30\penalty0 (7):\penalty0 854--867, 1984.
\newblock \doi{10.1287/mnsc.30.7.854}.

\bibitem[Pritsker et~al.(1969)Pritsker, Watters, and Wolfe]{pritsker1969multiproject}
A.~Alan~B. Pritsker, Lawrence~J. Watters, and Philip~M. Wolfe.
\newblock Multiproject scheduling with limited resources: A zero-one programming approach.
\newblock \emph{Management Science}, 16\penalty0 (1):\penalty0 93--108, 1969.
\newblock \doi{10.1287/mnsc.16.1.93}.

\bibitem[Proth and Sauer(1998)]{proth1998scheduling}
Jean-Marie Proth and Nathalie Sauer.
\newblock Scheduling of piecewise constant product flows: A petri net approach.
\newblock \emph{European Journal of Operational Research}, 106\penalty0 (1):\penalty0 45--56, 1998.
\newblock \doi{10.1016/S0377-2217(98)00273-2}.

\bibitem[Sprecher(2000)]{sprecher2000scheduling}
Arno Sprecher.
\newblock Scheduling resource-constrained projects competitively at modest memory requirements.
\newblock \emph{Management Science}, 46\penalty0 (5):\penalty0 710--723, 2000.
\newblock \doi{10.1287/mnsc.46.5.710.12044}.

\bibitem[van~der Aalst and Carmona(2022)]{van2022process}
Wil M.~P. van~der Aalst and Josep Carmona, editors.
\newblock \emph{Process Mining Handbook}.
\newblock Springer, Cham, 2022.
\newblock \doi{10.1007/978-3-031-08848-3}.

\bibitem[Yousefzadeh and Harati~Motlagh(2026)]{yo26mmr}
H.~R. Yousefzadeh and S.~Harati~Motlagh.
\newblock Multi-mode resource constrained project scheduling in home health care: a {DEA} framework.
\newblock \emph{Journal of the Operational Research Society}, 77\penalty0 (8):\penalty0 2311--2329, 2026.
\newblock \doi{10.1080/01605682.2025.2564768}.

\bibitem[Zhang et~al.(2024)Zhang, Li, Demeulemeester, and Zhang]{zh24atd}
J.~Zhang, L.~Li, E.~Demeulemeester, and H.~Zhang.
\newblock A three-dimensional spatial resource-constrained project scheduling problem: Model and heuristic.
\newblock \emph{European Journal of Operational Research}, 319\penalty0 (3):\penalty0 943--966, 2024.
\newblock \doi{10.1016/j.ejor.2024.07.018}.

\bibitem[Čapek et~al.(2012)Čapek, Šůcha, and Hanzálek]{ca12psw}
R.~Čapek, P.~Šůcha, and Z.~Hanzálek.
\newblock Production scheduling with alternative process plans.
\newblock \emph{European Journal of Operational Research}, 217\penalty0 (2):\penalty0 300--311, 2012.
\newblock \doi{10.1016/j.ejor.2011.09.018}.

\end{thebibliography}

\clearpage
\appendix
\numberwithin{equation}{section}
\counterwithin{figure}{section}
\counterwithin{table}{section}
\section*{Appendices}
\section{Proofs}\label{app:proofs}
In this appendix, we present the proofs of the propositions in Section~\ref{sec:equivalences}.

\subsection{Proof of Proposition~\texorpdfstring{\ref{prop:sa_pn}}{1}}\label{app:proof_prop1}

\begin{restatedprop}{\ref{prop:sa_pn}}
The state graph $G=(\Sset,\mathbb{E})$ and the $TTPNR$ reachability graph
$RG=(\Sset_{RG},\mathbb{E}_{RG})$, where $\Sset_{RG}$ is the set of reachable timed states and $\mathbb{E}_{RG}$ is the set of directed edges induced by transition firings, are isomorphic. That is, there exists a one-to-one correspondence between the states of $G$ and the states of $RG$, $\phi:\Sset\to\Sset_{RG}$, such that:
\begin{enumerate}
    \item[\emph{(i)}] $\phi(s_0)$ is the initial RG node represented by the timed marking $M_{s_0}$ and system time $0$, where $M_{s_0}$ refines the ordinary initial marking $M_0$,
    \item[\emph{(ii)}] corresponding states have the same system time,
    \item[\emph{(iii)}] for every activity $j\in\V$, there is an edge $s \xrightarrow{\,j\,} s'$ in $G$ if and only if there is an edge $\phi(s) \xrightarrow{\,t_j\,} \phi(s')$ in $RG$ with label $L(t_j)=j$.
\end{enumerate}
\end{restatedprop}

\begin{proof}
\red{For each $s\in\Sset$, we define $\phi(s)=(g(s),M_s)$, where $M_s$ is constructed below.}


For each control-flow place \red{$p\in
P_C\setminus\{p_{\mathrm{start}},p_{\mathrm{end}}\}$}, let $\bullet p=\{t_i\}$ and $p\bullet=\{t_j\}$ denote its unique input and output transitions. We set

\begin{equation}
M_s(p)=
\begin{cases}
\{1^{\,\dres{i}(s)}\},
& \text{if } i\in\V^o(s),\\[-.2cm]
\{1^0\},
& \text{if } i\in\V^c(s)
  \text{ and } j\notin\V^c(s)\cup\V^o(s),\\[-.2cm]
\emptyset,
& \text{otherwise.}
\end{cases}
\end{equation}

\red{For the source place, $M_s(p_{\mathrm{start}})=\{1^0\}$ if $0\notin\V^c(s)$, and
$\emptyset$ otherwise. For the sink place, $M_s(p_{\mathrm{end}})=\{1^0\}$ if $n\in\V^c(s)$, and $\emptyset$ otherwise.}

Thus, if activity $i$ is ongoing, then its output token toward each successor is present as a delayed token with residual duration $\dres{i}(s)$. If activity $i$ has completed and successor $j$ has not yet started, then the corresponding place contains one immediately available token. Otherwise the place is empty. 

For each resource place $p_{r_k}\in P_R$, the timed multiset $M_s(p_{r_k})$ consists of $R_k-\sum_{i\in\V^o(s)} r_{ik}$ tokens with residual duration $0$ and, for each ongoing activity $i\in\V^o(s)$ with $r_{ik}>0$, exactly $r_{ik}$ tokens with residual duration $\dres{i}(s)$. Hence $M_s(p_{r_k})$ records both the resource units that are immediately available at time $g(s)$ and the resource units that are currently tied to ongoing activities and return after their remaining processing times.

The construction above determines $\phi(s)$ uniquely from the state tuple 
$\bigl(g(s),\,\V^c(s),\,\V^o(s),\,\tilde d(s)\bigr)$.
By construction, $\phi(s)$ has the same system time as $s$, so part \emph{(ii)} holds. Moreover, for the initial state $s_0$ there are no completed or ongoing activities and all resources are available, hence $\phi(s_0)$ is the initial RG node represented by $M_{s_0}$ and system time $0$, where $M_{s_0}$ is obtained from $M_0$ by assigning zero residual availability time to every token, proving part \emph{(i)}.

Conversely, for the RCPSP-induced nets considered here, every reachable timed state determines the corresponding state-action information uniquely. Delayed control-flow tokens identify activities that have started but whose output places are not yet available, and their residual times determine the
corresponding values $\dres{i}(s)$. Immediately available control-flow tokens indicate completed predecessors whose successor activities have not yet consumed the corresponding tokens. The resource-place tokens distinguish between immediately available resource units and resource units committed to ongoing activities. \red{ Starting from activity $0$, whose status is determined by the source-place token, we recover the remaining activity statuses in topological order where an activity has started if and only if all its predecessors have completed and all its input control-flow places are empty. A started activity that is not ongoing is completed.} Since only reachable states are considered, these token patterns arise from a unique combination of completed activities, ongoing activities, residual durations, and system time. Thus the correspondence is one-to-one on the reachable state space.

It remains to prove part \emph{(iii)}. \red{By the construction of $M_s$, the control-flow tokens required by $t_j$, possibly delayed, are present if and only if $j\in\V^a(s)$. Resource constraints affect only the earliest firing time. Thus, the selectable transitions correspond exactly to $\V^a(s)$.}

In the state graph, the successor $s'$ is obtained by advancing the clock by the smallest amount $\Delta = g(s')-g(s)\ge 0$ such that all predecessors of $j$ have completed and enough resources are available to start $j$.

Now consider the RG node $\phi(s)$, represented by the timed marking
$M_s$ and current system time $g(s)$. For every predecessor $i$ of $j$, the corresponding control-flow place on an input arc of $t_j$ contains exactly one token: it is immediately available when $i\in\V^c(s)$, and it is delayed by $\dres{i}(s)$ when $i\in\V^o(s)$. Likewise, each resource place contains exactly the immediately available resource tokens together with the delayed return tokens induced by the ongoing activities in $s$. Therefore, the smallest time increment after which all input control-flow and resource tokens required by $t_j$ are available according to $M_s$ is exactly the same value $\Delta$ as in the state-action transition rule.

After advancing time by $\Delta$ and firing $t_j$, the Petri-net token dynamics match the state-action update exactly:
\begin{itemize}
    \item the required predecessor and resource tokens for $j$ are consumed,
    \item if $d_j>0$, activity $j$ becomes ongoing, so its output control-flow tokens and its used resource tokens reappear as delayed tokens with residual duration $d_j$; if $d_j=0$, activity $j$ completes immediately and its output tokens and resources are immediately available,
    \item for every previously ongoing activity, residual durations are reduced by $\Delta$ (but not below 0),
    \item activities whose residual duration becomes zero leave $\V^o$ and join $\V^c$, and their successor places now contain immediately available tokens unless these have already been consumed by already started successors.
\end{itemize}
Hence the resulting timed marking is precisely $M_{s'}$, and the resulting RG node is $\phi(s')$. This proves that every edge $s \xrightarrow{\,j\,} s'$ in $G$ corresponds to an edge $\phi(s) \xrightarrow{\,t_j\,} \phi(s')$ in $RG$.

The converse follows by the same argument. Every edge in $RG$ is defined by choosing one transition $t_j$, advancing time by the smallest amount needed to make its input control-flow and resource tokens available, and then firing it. This is exactly the state-transition rule used in Section~\ref{sec:dp} for the scheduling decision $j$. Therefore each edge in $RG$ corresponds to a unique edge in $G$ with the same label.
Thus, $\phi$ preserves and reflects labeled adjacency. 
Since both graphs start from corresponding initial states and generate exactly the same labeled successors from corresponding states, \red{induction on path length establishes that
$\phi(\Sset)=\Sset_{RG}$.}
\end{proof}

\subsection{Proof of Proposition~\texorpdfstring{\ref{prop:milp_sa}}{2}}\label{app:proof_prop2}

\begin{restatedprop}{\ref{prop:milp_sa}}
    Given an RCPSP, let $G$ be the state graph associated with the state-action representation. Then, any partial schedule $(\J,\{S_i\}_{i\in\J})$ produced by earliest-start generation corresponds to a path in the graph $G$ leading to a state $s\in\Sset$ satisfying
    $g(s)=\max_{i\in\J} S_i$ and $\V^o(s)=\{i\in\J: S_i\leq g(s)< S_i+d_i\}$, $\tilde{d}_i(s)=S_i+d_i-g(s)>0$ for all $i\in \V^o(s)$, and $\V^c(s)=\J\setminus \V^o(s)$.
    Moreover, any state $s\in \Sset$ corresponds to at least one partial schedule $(\J,\{S_i\}_{i\in\J})$ such that $\J=\V^o(s)\cup\V^c(s)$, and $S_i=g(s)+\tilde{d}_i(s)-d_i$ for all $i\in \V^o(s)$ and $S_i\leq g(s)-d_i$ for all $i\in \V^c(s)$.
\end{restatedprop}
\begin{proof}
\red{Since $(\J,S_\J)$ is produced by earliest-start generation, there exists a sequence $(i_0=0,i_1,\ldots,i_k)$ of the activities in $\J$ such that, for each $\ell\ge1$, $S_{i_\ell}$ is the earliest start time of $i_\ell$ that is no earlier than $S_{i_{\ell-1}}$ and feasible with respect to $(\{i_0,\ldots,i_{\ell-1}\},\{S_{i_0},\ldots,S_{i_{\ell-1}}\})$, i.e., $S_{i_\ell}$ is given by \eqref{eq:compute_earlier_start}.
We denote this ordered list as $\J^o=(i_0=0,i_1,\ldots,i_k)$.
We can therefore construct the path in $G$ starting from $s_0$ by applying the scheduling decisions from $\J^o$ in order. 
Specifically, at state $s_{\ell}$ we apply decision $i_{\ell}$ and proceed to state $s_{\ell+1}$ by applying the state transition formulas. 
By induction on $\ell\geq 1$, using \eqref{eq:ongoing_transition} and \eqref{eq:completed_transition}, each state $s_\ell$ satisfies $g(s_\ell)=S_{i_{\ell-1}}$, and $\V^o(s_\ell)$ consists of the activities among $i_0,\ldots,i_{\ell-1}$ that are in process at $g(s_\ell)$, with $\tilde d_i(s_\ell)=S_i+d_i-g(s_\ell)>0$, and $\V^c(s_\ell)=\{i_0,\ldots,i_{\ell-1}\}\setminus\V^o(s_\ell)$. Hence the conditions in \eqref{eq:time_transition} coincide with those in \eqref{eq:compute_earlier_start}, and $g(s_{\ell+1})=S_{i_\ell}$.}
Moreover, we terminate the path in $s_{k+1}$, which must satisfy $g(s_{k+1})=S_{i_k}=\max_{j\in\J} S_j$. The definitions of $\V^o(s_{k+1})$ and $\V^c(s_{k+1})$ follow from the recursive transition.

Conversely, let $s\in \Sset$, and let $(s_0,s_1,\ldots,s_{k}=s)$ be a path from $s_0$ to $s$ on graph $G$, then by definition there exist scheduling decisions $(j_0=0,j_1,j_2,\ldots, j_{k-1})$ such that $j_\ell$ induces the transition from $s_{\ell}$ to $s_{\ell+1}$.
Moreover, since $j_\ell$ is feasible for state $s_\ell$ the set $\{j_0,j_1,\ldots,j_\ell\}$ is a predecessor-closed set, and together with $S_{j_\ell}=g(s_{\ell+1})$ it defines a partial schedule.
\red{Since each $S_{j_\ell}$ is the earliest feasible start time given $(\{j_0,\ldots,j_{\ell-1}\},\{S_{j_0},\ldots,S_{j_{\ell-1}}\})$, this partial schedule is produced by earliest-start generation, and the first part of the proof yields the stated relations between $(\J,S_\J)$ and $s$.}
\end{proof}

\section{Empirical study}\label{app:empirical}
In this appendix, we present the LP relaxations used in the empirical study of Section~\ref{sec:empirical} as well as additional computational results.

\subsection{Weak and strong LP relaxations}\label{app:lp}
Both LP relaxations are derived from the time-indexed MILP formulation of Section~\ref{sec:milp}, whose notation we adopt throughout. In both relaxations, the integrality requirements are relaxed, i.e., $y_{it}\in\{0,1\}$ is replaced by $0\le y_{it}\le 1$ and $S_i\in\mathbb{N}_0$ by $S_i\ge 0$, while the objective, the assignment constraints, and the start-time definitions are kept unchanged. The two relaxations differ only in how the precedence relations are modeled.

\paragraph{Node-specific LP and resulting bounds}
Consider a B\&B node $N_i$ with partial schedule $(\J,S_{\J})$ and current time $g(N_i)$. Instead of the full horizon $\T=\{0,\dots,H\}$, the LP at $N_i$ only considers the time window 
\begin{equation}
  \T(N_i)=\{g(N_i),\dots,\bar F\},
\end{equation}
which runs from the current time $g(N_i)$, a lower bound on the start time of every unscheduled activity, to the global upper bound $\bar F$, i.e., the makespan of the incumbent. The upper bound $\bar F$ is initialized to the horizon $H$ and is updated whenever a new incumbent is found. Since we are only interested in schedules that start no unscheduled activity before $g(N_i)$ and have makespan at most $\bar F$, this restriction does not cut off any schedule in the subtree of $N_i$ that can match or improve the incumbent. Accordingly, the LP at $N_i$ has decision variables $y_{jt}$ only for unscheduled activities $j\in\V\setminus\J$ and start times $t\in\{g(N_i),\dots,\bar F-d_j\}$. All other values of $y_{jt}$ are set to constants:
\begin{itemize}
  \item for scheduled activities $j\in\J$, whose start times $S_j$ are given, $y_{jt}=1$ if $t=S_j$ and $y_{jt}=0$ otherwise;
  \item for unscheduled activities $j\in\V\setminus\J$, $y_{jt}=0$ for $t<g(N_i)$ and for $t>\bar F-d_j$.
\end{itemize}
Substituting these values into the LP affects the constraints as follows.
\begin{itemize}
  \item \emph{Assignment constraints} are only needed for unscheduled activities, and read $\sum_{t=g(N_i)}^{\bar F-d_j} y_{jt}=1$ for all $j\in\V\setminus\J$.
  \item \emph{Resource constraints} are imposed only for $t\in\T(N_i)$, with the capacity reduced by the scheduled activities still in process at time $t$, i.e., the right-hand side becomes $R_k-\sum_{j\in\J:\,S_j\le t<S_j+d_j} r_{jk}$. For $t<g(N_i)$, only scheduled activities are in process, so these constraints are satisfied by the feasibility of the partial schedule and can be omitted.
  \item \emph{Precedence constraints} $(i,j)\in\E$ with $i,j\in\J$ involve only constants, are satisfied by the partial schedule, and can be omitted. Since an activity is only scheduled after all of its predecessors have been scheduled, the case $i\notin\J$, $j\in\J$ cannot occur. Hence, only two cases remain: precedence relations between two unscheduled activities, $i,j\notin\J$, and precedence relations in which only the predecessor is scheduled, $i\in\J$ and $j\notin\J$. In the latter case, the constraint reduces to a release time for $j$. The resulting constraints are given below for each relaxation.
\end{itemize}
The purpose of the time window is to reduce the size of the LP formulation, as the number of variables, resource constraints, and (for the strong LP) precedence constraints all grow linearly with the length of the window, so a shorter window yields a smaller LP that is faster to solve.

Let $z_{\mathrm{LP}}(N_i)$ denote the optimal value of the LP at node $N_i$ (weak or strong, see below). Since every start time in the window is at least $g(N_i)$, we have $z_{\mathrm{LP}}(N_i)\ge g(N_i)$, so the LP provides a lower bound $z_{\mathrm{LP}}(N_i)-g(N_i)$ on the remaining time to project completion. Because makespans are integral, this bound can be rounded up and combined with the critical-path bound $h(N_i)$ of Section~\ref{sec:approaches}, yielding
\begin{equation}
  h_{\mathrm{LP}}(N_i)=\max\Bigl\{h(N_i),\
  \bigl\lceil z_{\mathrm{LP}}(N_i)\bigr\rceil-g(N_i)\Bigr\}.
\end{equation}
In the LP-based variants, the evaluation function of Section~\ref{sec:approaches} is replaced by $\feval(N_i)=g(N_i)+h_{\mathrm{LP}}(N_i)$, which is used both for node selection and for pruning. A node whose LP is infeasible is pruned as well. Moreover, if the optimal LP solution at $N_i$ is integral, i.e., $y_{jt}\in\{0,1\}$ for all $j\in\V\setminus\J$ and $t\in\T(N_i)$, it satisfies all constraints of the MILP formulation and, together with the partial schedule $(\J,S_{\J})$, defines a feasible schedule with makespan $z_{\mathrm{LP}}(N_i)$. In this case, the upper bound is updated to $\bar F\leftarrow\min\{\bar F,z_{\mathrm{LP}}(N_i)\}$, the incumbent is updated accordingly, and $N_i$ is pruned, since no schedule in its subtree can have a smaller makespan.

\paragraph{Weak LP relaxation}
The weak LP keeps the precedence constraints of the MILP formulation as they are, i.e., each precedence relation $(i,j)\in \E$ yields the single aggregated constraint
\begin{equation}
  S_i+d_i\le S_j
  \quad\Longleftrightarrow\quad
  \sum_{t=0}^{\bar F-d_i} (t+d_i)\,y_{it} \;\le\; \sum_{t=0}^{\bar F-d_j} t\,y_{jt}
  \qquad \forall (i,j)\in \E. \label{eq:lp-weak}
\end{equation}
Hence, the weak LP is the plain LP relaxation of the MILP formulation.

\paragraph{Strong LP relaxation}
The strong LP replaces~\eqref{eq:lp-weak} by one cumulative constraint per precedence relation and time point $t\in\T(N_i)$: by any time $t$, activity $j$ may have started only to the extent that activity $i$ has already completed. If $i$ is already scheduled, this simply prohibits $j$ from starting before $i$ completes. Hence, \eqref{eq:lp-weak} is replaced by
\begin{align}
  \sum_{\tau=g(N_i)}^{t-d_i} y_{i\tau} &\;\ge\;
  \sum_{\tau=g(N_i)}^{\min\{t,\bar F-d_j\}} y_{j\tau}
  && \forall (i,j)\in \E,\ i,j\notin\J,\ t\in\T(N_i), \label{eq:lp-strong}\\
  \sum_{t=g(N_i)}^{S_i+d_i-1} y_{jt} &\;\le\; 0
  && \forall (i,j)\in \E,\ i\in\J,\ j\notin\J. \label{eq:lp-strong-release}
\end{align}
Together with the assignment constraints, summing~\eqref{eq:lp-strong} over $t\in\T(N_i)$ yields~\eqref{eq:lp-weak} for $i,j\notin\J$, and~\eqref{eq:lp-strong-release} implies~\eqref{eq:lp-weak} for $i\in\J$,
$j\notin\J$, since it allows $j$ to start only at or after the completion time $S_i+d_i$ of $i$. Hence, every solution of the strong LP satisfies the weak one and the strong bound is at least as large. The price is
up to $|\E|\,|\T(N_i)|$ precedence constraints instead of $|\E|$.

\subsection{Additional computational results}\label{app:additional}
\begin{figure}[t]
    \centering
    \includegraphics[width=\linewidth]{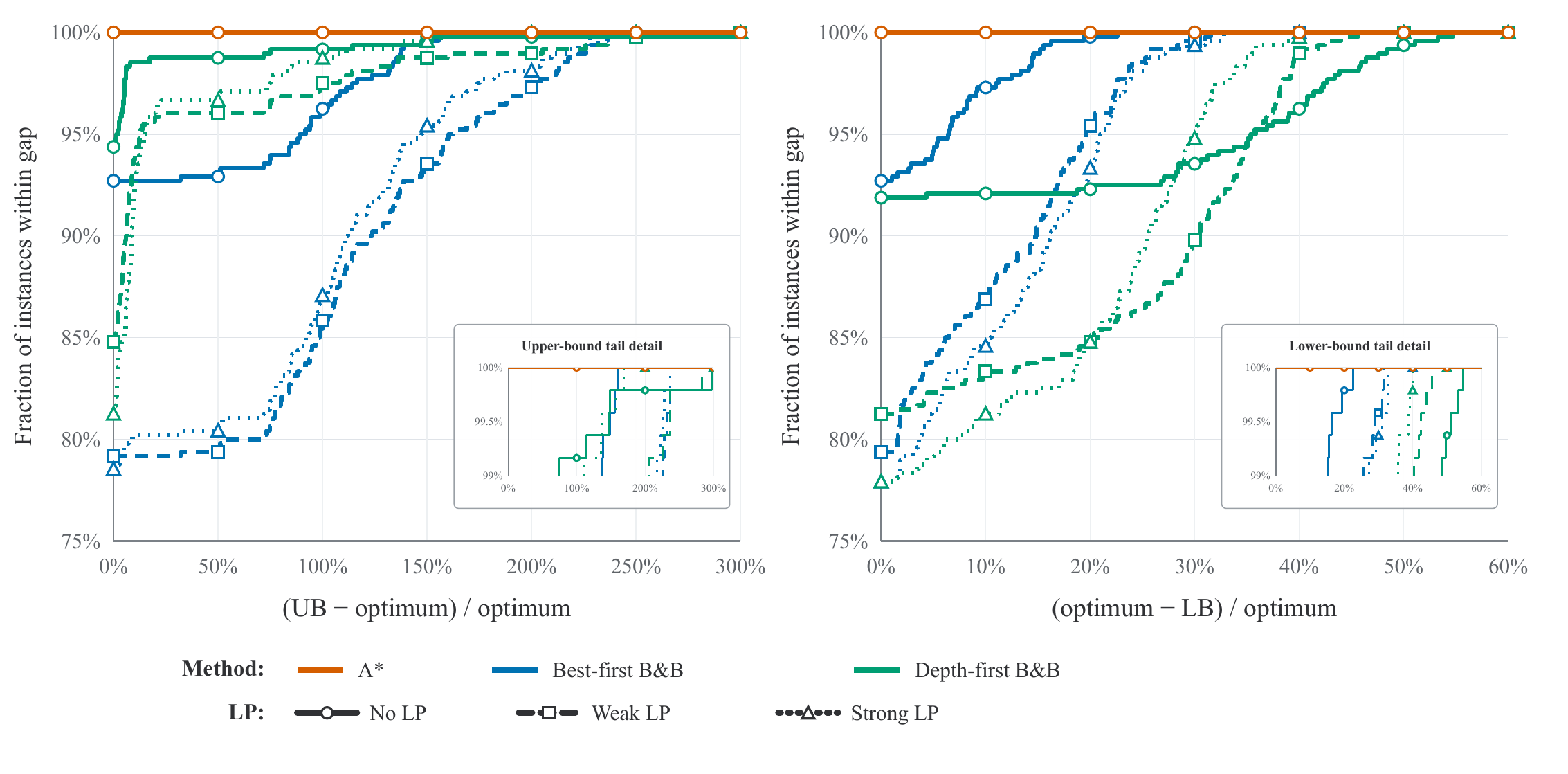}
    \caption{Profiles of the final upper bounds (left) and lower bounds (right) over all 480 instances, }
    \label{fig:bb-performance-profile}
    \vspace{-.1cm}\caption*{\scriptsize \textit{Note:} Bounds measured as the relative gaps (UB$-$optimum)/optimum and (optimum$-$LB)/optimum, respectively, in percent. Each curve gives the fraction of instances whose gap does not
    exceed the given value. Since $A^*$ proves optimality on all instances, its curves are at one for every gap.}
\end{figure}

\begin{figure}[t]
    \centering
    \includegraphics[width=\linewidth]{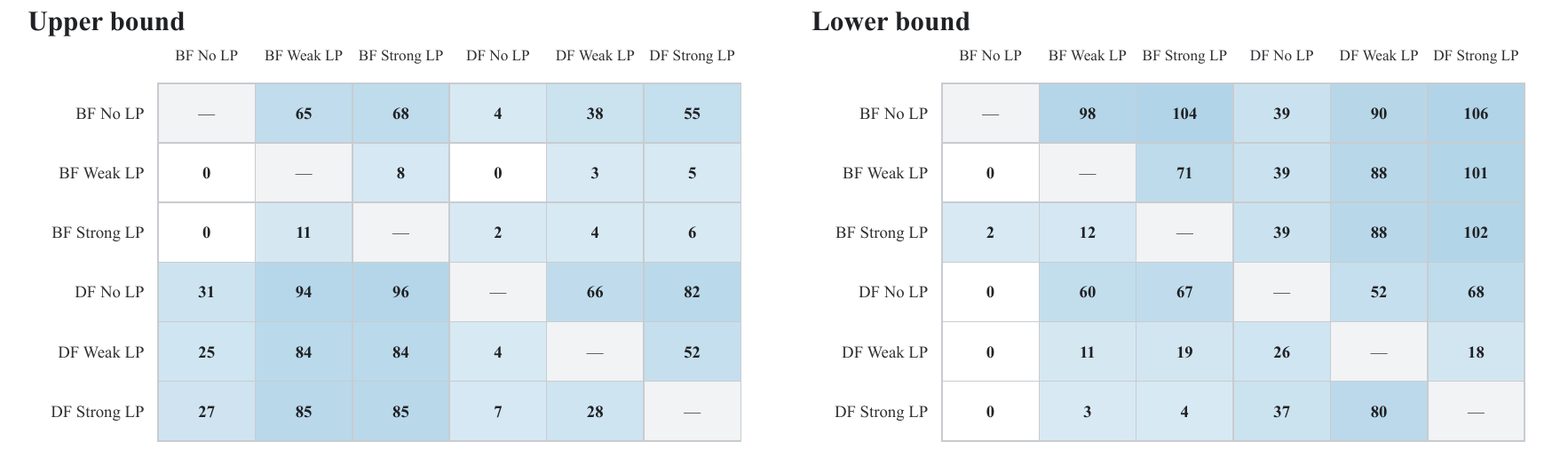}
    \caption{Pairwise comparison of the final upper bounds (left) and lower bounds (right) obtained by the B\&B variants.}
    \label{fig:pairwise-win-tie-loss}
    \vspace{-.1cm}\caption*{\scriptsize \textit{Note:} $A^*$ is not included, as it proves optimality on all instances. Each cell gives the number of instances on which the row method obtained a better bound than the column method; darker cells indicate that the row method dominates the column method on more instances.}
\end{figure}

This section complements the empirical study of Section~\ref{sec:empirical}. Comparing the B\&B variants among themselves is not the main goal of the paper, which focuses on the comparison between B\&B and $A^*$. The following results are therefore reported for the interested reader. Figure~\ref{fig:bb-performance-profile} compares the final upper and lower bounds obtained by each method with the optimal makespan (for $A^*$, both
bounds coincide with the optimum on all instances), and Figure~\ref{fig:pairwise-win-tie-loss} reports, for each pair of B\&B variants, on how many instances one variant obtained a better upper or lower
bound than the other. Since $A^*$ proves optimality on all instances, it is not included in the pairwise comparison. Here, UB denotes the final upper bound $\bar F$, i.e., the makespan of the best schedule found, and LB the final lower bound, i.e., the smallest lower bound $\feval(N_i)$ among the nodes that remain open when the search terminates.

Both figures show that DF B\&B is better at obtaining upper bounds, whereas BF B\&B is better at obtaining lower bounds. For instance, without LP relaxations, DF B\&B obtains a better upper bound than BF B\&B on 31
instances and a worse one on only 4, whereas BF B\&B obtains a better lower bound on 39 instances and never a worse one. The LP relaxations generally deteriorate the upper bounds, since their computational overhead slows down the search and fewer incumbents are found within the time limit. For example, DF B\&B without LP obtains a better upper bound than with the weak and strong LP on 66 and 82 instances, respectively, and a worse one on only 4 and 7 instances. The LP relaxations also do not improve the lower bounds of BF B\&B, as the no-LP variant obtains a better lower bound than both LP variants on at least 98 instances. For DF B\&B, in contrast, the LP relaxations present a trade-off, which is why their lower-bound profiles cross that of the no-LP variant. The no-LP variant proves optimality on more instances, since solving an LP at every node is two to three orders of magnitude slower, but the LP variants obtain considerably tighter lower bounds on the instances that remain unsolved. This trade-off is sharper for the strong LP, which proves optimality on 69 fewer instances than the no-LP variant, compared with 52 for the weak LP, but also yields tighter lower bounds. Specifically, on the 40 instances that the no-LP variant does not solve, which the LP variants do not solve either, the median relative gap between the optimum and the final lower bound drops from 39.0\% without LP relaxations to 34.2\% with the weak LP and 28.1\% with the strong LP. On these instances, the final lower bound of DF B\&B usually remains at the root bound, as shallow open nodes are never revisited. The strong LP therefore benefits mainly from its tighter root bound, whereas the weak LP also raises its bound above the root bound on about a third of these instances, as its cheaper LPs allow more nodes to be processed. This suggests that, for DF B\&B, it may be preferable to solve the strong LP only once at the root, to obtain a global lower bound, and to use the cheaper critical-path bound at all other nodes. Such a hybrid would retain most of the lower-bound advantage of the strong LP while preserving the speed of the no-LP variant.

\begin{table*}[t]
  \centering
  \caption{Instances for which optimality was not proven, partitioned by PSPLIB parameter subclass.}
  \label{tab:unresolved-by-psplib-subclass}
  \footnotesize
  \renewcommand{\arraystretch}{0.92}
  \begin{tabular}{@{}rrr r rrr rrr@{}}
    \toprule
    & & & & \multicolumn{3}{c}{Best-first B\&B}
          & \multicolumn{3}{c}{Depth-first B\&B} \\
    \cmidrule(lr){5-7}\cmidrule(l){8-10}
    $RS$ & $RF$ & $NC$ & A$^*$
      & No LP & Weak LP & Strong LP
      & No LP & Weak LP & Strong LP \\
    \midrule
    0.25 & 0.25 & 1.5 & 0 &  0 &  0 &  0 &  0 &  0 &  1 \\
    0.25 & 0.25 & 1.8 & 0 &  0 &  0 &  1 &  0 &  0 &  1 \\
    0.25 & 0.50 & 1.5 & 0 &  1 &  9 &  9 &  1 &  8 &  9 \\
    0.25 & 0.50 & 1.8 & 0 &  0 &  5 &  9 &  0 &  3 &  7 \\
    0.25 & 0.50 & 2.1 & 0 &  0 &  4 &  6 &  0 &  4 &  6 \\
    0.25 & 0.75 & 1.5 & 0 &  5 & 10 & 10 &  6 & 10 & 10 \\
    0.25 & 0.75 & 1.8 & 0 &  5 & 10 & 10 &  5 & 10 & 10 \\
    0.25 & 0.75 & 2.1 & 0 &  2 & 10 & 10 &  3 & 10 & 10 \\
    0.25 & 1.00 & 1.5 & 0 &  9 & 10 & 10 & 10 & 10 & 10 \\
    0.25 & 1.00 & 1.8 & 0 &  9 & 10 & 10 &  9 & 10 & 10 \\
    0.25 & 1.00 & 2.1 & 0 &  4 & 10 & 10 &  4 & 10 & 10 \\
    \addlinespace
    0.50 & 0.50 & 1.8 & 0 &  0 &  1 &  1 &  0 &  0 &  1 \\
    0.50 & 0.75 & 1.5 & 0 &  0 &  4 &  4 &  1 &  4 &  3 \\
    0.50 & 0.75 & 1.8 & 0 &  0 &  0 &  0 &  0 &  1 &  1 \\
    0.50 & 0.75 & 2.1 & 0 &  0 &  1 &  1 &  0 &  0 &  1 \\
    0.50 & 1.00 & 1.5 & 0 &  0 &  5 &  6 &  1 &  6 &  6 \\
    0.50 & 1.00 & 1.8 & 0 &  0 &  5 &  5 &  0 &  3 &  4 \\
    0.50 & 1.00 & 2.1 & 0 &  0 &  3 &  4 &  0 &  1 &  4 \\
    \addlinespace
    0.75 & 0.75 & 1.5 & 0 &  0 &  0 &  0 &  0 &  0 &  1 \\
    0.75 & 1.00 & 1.5 & 0 &  0 &  2 &  2 &  0 &  2 &  2 \\
    0.75 & 1.00 & 1.8 & 0 &  0 &  1 &  1 &  0 &  0 &  2 \\
    \bottomrule
  \end{tabular}%
  \caption*{\scriptsize \textit{Note:} Each $(NC,RF,RS)$ subclass contains ten instances, so every entry is a count out of ten. Subclasses omitted from the table have zero unresolved instances for every method.}
\end{table*}

Table~\ref{tab:unresolved-by-psplib-subclass} reports the instances for which optimality was not proven, partitioned by PSPLIB parameter subclass. The J30 instances are generated by varying three parameters \citep{kolisch1997psplib}. The network complexity ($NC$) is the average number of immediate successors per activity. The resource factor ($RF$) is the average fraction of the resources required by an activity, so that a higher $RF$ means that each activity uses more of the resources. The resource strength ($RS$) measures how scarce the resources are relative to the demand, where a lower $RS$ means more scarce resources and hence a more constrained instance. The unresolved instances concentrate in subclasses with low resource strength ($RS=0.25$) and high resource factor ($RF\ge 0.75$), and the LP-based variants leave considerably more instances unresolved than the variants without LP relaxations.

\end{document}